\documentclass[reqno,12pt]{amsart}

\usepackage{color} 
\usepackage{ifpdf}
\ifpdf
    \usepackage[pdftex]{graphicx}
    \usepackage[pdftex]{hyperref}
    \hypersetup{
        unicode=false,          
        pdftoolbar=true,        
        pdfmenubar=true,        
        pdffitwindow=false,     
        pdfstartview={FitH},    
        pdftitle={MCP Article},      
        pdfauthor={Michael Holst},   
        pdfsubject={Mathematics},    
        pdfcreator={Michael Holst},  
        pdfproducer={Michael Holst}, 
        pdfkeywords={PDE, analysis, mathematical physics}, 
        pdfnewwindow=true,      
        colorlinks=true,        
        linkcolor=red,          
        citecolor=blue,         
        filecolor=magenta,      
        urlcolor=cyan           
    }

\else
    \usepackage{graphicx}
    \usepackage{pstricks}
    
    \newcommand{\href}[2]{#2}
\fi

\usepackage{times}
\usepackage{amsfonts}

\usepackage{amsmath}
\usepackage{amsthm}
\usepackage{amssymb}
\usepackage{amsbsy}
\usepackage{amscd}

\usepackage{enumerate}
\usepackage{verbatim}
\usepackage{subfigure}

\newtheorem{theorem}{Theorem}[section]

\newtheorem{lemma}[theorem]{Lemma}

\newtheorem{proposition}[theorem]{Proposition}

\newtheorem{remark}[theorem]{Remark}

\numberwithin{equation}{section}  

  \newcounter{mnote}
  \let\oldmarginpar\marginpar
    \renewcommand\marginpar[1]{\-\oldmarginpar[\raggedleft\footnotesize #1]%
    {\raggedright\footnotesize #1}}

\definecolor{myblue}{rgb}{0.2,0.2,0.7}
\definecolor{mygreen}{rgb}{0,0.6,0}
\definecolor{mycyan}{rgb}{0,0.6,0.6}
\definecolor{myred}{rgb}{0.9,0.2,0.2}
\definecolor{mymagenta}{rgb}{0.9,0.2,0.9}
\definecolor{mywhite}{rgb}{1.0,1.0,1.0}
\definecolor{myblack}{rgb}{0.0,0.0,0.0}

\usepackage[curve]{xypic}

\newcommand{\poly} {{\mathcal P}}

\newcommand{\R} {\mathbb R}

\newcommand{\Div}{\text{div}}
\newcommand{\Curl}{\text{curl}}
\newcommand{\p}{\partial}

\newcommand{\cancel}[1]{}
\newcommand{\raw}{\rightarrow}

\newcommand{\Tau}{\mathcal T}
\newcommand{\eps}{\epsilon}
\newcommand{\veps}{\varepsilon}

\newcommand{\ol}[1]{\overline{#1}}
\newcommand{\mr}[1]{\mathring{#1}}

\newcommand{\id}{\mathbb{I}}

\newcommand{\vn}[1]{\left|\left|#1\right|\right|}

\newcommand{\ds}{\displaystyle}

\newcommand{\eltn}[1]{{\vn{#1}_{L^2}}}

\newcommand{\hdiv}{H(\textnormal{div})}

\newcommand{\ut}{\mu}

\newcommand{\nedelec}{N\'ed\'elec }

\begin{document}

\title[FEEC for Evolution Problems]   
      {Finite Element Exterior Calculus for Evolution Problems}

\author[A. Gillette]{Andrew Gillette}
\email{agillette@math.arizona.edu}
\address{Department of Mathematics\\
         University of Arizona\\ 
         Tucson AZ 85721}
         
\author[M. Holst]{Michael Holst}
\email{mholst@math.ucsd.edu}
\address{Department of Mathematics\\
         University of California San Diego\\ 
         La Jolla CA 92093}

\author[Y. Zhu]{Yunrong Zhu}
\email{zhuyunr@isu.edu}
\address{Department of Mathematics \& Statistics\\
         Idaho State University\\
         Pocatello ID 83209}

\date{\today}


\begin{abstract}
  Arnold, Falk, and Winther
  [\emph{Bull.\ Amer.\ Math.\ Soc.}~\textbf{47} (2010), 281--354]
  showed that mixed variational problems, and their
  numerical approximation by mixed methods, could be most completely 
  understood using the ideas and tools of {\em Hilbert complexes}. 
  This led to the development of the Finite Element Exterior Calculus (FEEC) 
  for a large class of linear elliptic problems.
  More recently, Holst and Stern 
  [\emph{Found.\ Comp.\ Math.}~\textbf{12}:3 (2012), 263--293 and 363--387]
  extended the FEEC framework to semi-linear problems, and to problems 
  containing {\em variational crimes}, allowing for the analysis and numerical 
  approximation of linear and nonlinear geometric elliptic partial 
  differential equations on Riemannian manifolds of arbitrary spatial 
  dimension, generalizing surface finite element approximation theory.
  In this article, we develop another distinct extension to the FEEC, namely
  to parabolic and hyperbolic evolution systems, allowing for the treatment
  of geometric and other evolution problems.
  Our approach is to combine the recent work on the FEEC for elliptic problems 
  with a classical approach to solving evolution problems via semi-discrete
  finite element methods, by viewing solutions to the evolution problem
  as lying in time-parameterized Hilbert spaces (or {\em Bochner} spaces).
  Building on classical approaches by Thom\'ee for parabolic problems 
  and Geveci for hyperbolic problems, we establish  
  {\em a priori} error estimates for Galerkin FEM approximation in the 
  natural parametrized Hilbert space norms.    
  In particular, we recover the results of Thom\'ee and Geveci for 
  two-dimensional domains and lowest-order mixed methods as special cases,
  effectively extending their results to arbitrary spatial dimension and to
  an entire family of mixed methods.
  We also show how the Holst and Stern framework allows for extensions of 
  these results to certain semi-linear evolution problems.
\end{abstract}

\maketitle


\tableofcontents

\clearpage

\section{Introduction}

More than two decades of research on linear mixed variational problems, and
their numerical approximation by mixed methods, recently culminated in the
seminal work of Arnold, Falk, and Winther~\cite{AFW2010}.
The authors show how such problems are most completely understood
using the ideas and tools of {\em Hilbert complexes}, leading to the
development of the Finite Element Exterior Calculus (FEEC) for elliptic
problems.
In two related articles~\cite{HS2010a,HS2010b}, Holst and Stern
extended the Arnold--Falk--Winther framework to include {\em variational
crimes}, allowing for the analysis and numerical approximation of linear 
and nonlinear geometric elliptic partial differential equations on 
Riemannian manifolds of arbitrary spatial dimension, generalizing the 
existing surface finite element approximation theory in several 
directions.

In this article, we extend the FEEC in another direction, namely
to parabolic and hyperbolic evolution systems.
Our approach is to combine the recent work on the FEEC for elliptic problems 
with a classical approach to solving evolution problems using semi-discrete
finite element methods, by viewing solutions to the evolution problem
as lying in time-parameterized Banach (or {\em Bochner}) spaces.
Building on classical approaches by Thom\'ee for parabolic problems 
and Geveci for hyperbolic problems, we establish  
{\em a priori} error estimates for Galerkin FEM approximation in the 
natural parametrized Hilbert space norms. 
In particular, we recover the results of Thom\'ee and Geveci for 
two-dimensional domains and the lowest-order mixed method as a special case,
effectively extending their results to arbitrary spatial dimension and to
an entire family of mixed methods.
We also show how the Holst and Stern framework allows for extensions of 
these results to certain semi-linear evolution problems.

To understand why the finite element exterior calculus (FEEC) has
emerged in a natural way to become a major mathematical tool in the
development of numerical methods for PDE, we recall one of the many examples
presented at length in~\cite{AFW2010}.
Consider the vector Laplacian over a domain $\Omega\subset\R^3$:
\[
- \Delta u := - \mbox{ grad }\mbox{ div } u + \mbox{ curl }\mbox{ curl } u,
\]
with the boundary conditions $u\cdot n =0$ and $\Curl~u\times n =0$ on $\partial \Omega$.
Given an $L^{2}$ vector field data $f$, the natural variational formulation of the problem is:
Find $u \in H(\mbox{curl};\Omega) \cap H_0(\mbox{div};\Omega)$ such that
\begin{equation}
\label{eq:non-mixed-vec-lap}
\hspace{-0.2cm}
\int_{\Omega} [(\nabla \cdot u) (\nabla \cdot v)
     + (\nabla \times u) \cdot (\nabla \times v)]~dx 
= \int_{\Omega} f \cdot v ~dx,
~~
\forall v \in H(\mbox{curl};\Omega) \cap H_0(\mbox{div};\Omega).
\end{equation}
By introducing an intermediate variable $\sigma := -\Div\; u$, a natural alternative formulation is the following {\em mixed} form:
find $(\sigma,u) \in H^1(\Omega) \times H(\mbox{curl};\Omega)$ such that
\begin{align}
\int_{\Omega} (\sigma \tau - u \cdot \nabla \tau) ~dx &= 0,
    \quad \forall \tau \in H^1(\Omega),
\label{eq:mixed-vec-lap-1}
\\
\int_{\Omega} [\nabla \sigma \cdot v 
    + (\nabla \times u) \cdot (\nabla \times v)] ~dx 
   &= \int_{\Omega} f \cdot v ~dx,
    \quad \forall v \in H(\mbox{curl};\Omega).
\label{eq:mixed-vec-lap-2}
\end{align}
Using the standard finite element approach based on the
non-mixed formulation (\ref{eq:non-mixed-vec-lap}) (e.g. using continuous piecewise linear vector functions) 
can yield incorrect results 
if the domain has certain geometric features (e.g.\ domains with re-entrant corners) 
 or topological features (e.g.\ non-simply connected domains).
A standard finite element approach based on the mixed formulation 
(\ref{eq:mixed-vec-lap-1})-(\ref{eq:mixed-vec-lap-2}), on the other hand,
suffers neither of these difficulties and typically works extremely well.  

The explanation for why one approach fails and the other succeeds lies in 
the fundamental mathematical structures underlying the finite element method. In case of the domain with re-entrant corners, the failure of the non-mixed approach is due to the non-density of $H^1(\Omega)$ in $H(\mbox{curl}) \cap H_0(\mbox{div})$, as discussed in~\cite[Remark 19]{Hiptmair2002} and~\cite[Section 2.3]{AFW2010}. This non-density leads to a problem of 
\textit{inconsistency}, i.e.\ that the discrete approximation of the operators
and data do not approximate the continuous problem correctly as
the mesh size is taken to zero. On the other hand, the
error due to topological features, such as non-simply connected domains, 
can be traced to the presence
of non-zero \textit{harmonic vector fields} on the domain, i.e.\ vector
fields that are both curl-free and divergence-free.
The mixed formulation turns out to be both consistent and respectful 
of non-zero harmonic vector fields while the standard formulation does not.
A natural question is then: What is an appropriate mathematical framework 
for understanding these problems abstractly so that
a methodical construction of ``good'' finite element methods
can be carried out for these and similar initial boundary value problems?

The answer turns out to be {\em Hilbert Complexes}.
Hilbert complexes were originally studied in \cite{BrLe1992} 
as a way to
generalize certain properties of elliptic complexes, particularly the
Hodge decomposition and other aspects of Hodge theory.  
  A \emph{Hilbert complex} $ \left( W , \mathrm{d} \right) $ consists
  of a sequence of Hilbert spaces $ W ^k $, along with closed,
  densely-defined linear maps $ \mathrm{d} ^k \colon V ^k \subset W ^k
  \rightarrow V ^{ k + 1 } \subset W ^{ k + 1 } $, possibly unbounded,
  such that $ \mathrm{d} ^k \circ \mathrm{d} ^{k-1} = 0 $ for each
  $k$.
  \begin{equation*}
    \xymatrix{  \cdots \ar[r] & V ^{ k - 1 } \ar[r]^-{ \mathrm{d} ^{k-1} }
      & V ^k \ar[r] ^-{ \mathrm{d} ^k } & V ^{ k + 1 } \ar[r] & \cdots }
  \end{equation*}
  This Hilbert complex is said to be \emph{bounded} if $ \mathrm{d} ^k
  $ is a bounded linear map from $ W ^k $ to $ W ^{ k + 1 } $ for each
  $k$, i.e., $ \left( W , \mathrm{d} \right) $ is a cochain complex in
  the category of Hilbert spaces.  It is said to be \emph{closed} if
  the image $ \mathrm{d} ^k V ^k $ is closed in $ W ^{ k + 1 } $ for
  each $k$.
It was shown in
\cite{AFW2006,AFW2010} that Hilbert complexes are also a convenient
abstract setting for mixed variational problems and their numerical
approximation by mixed finite element methods, providing the
foundation of a framework called \emph{finite element exterior
  calculus}.  This line of research is
the culmination of several decades of work on mixed finite element
methods and computational 
electromagnetics~\cite{B1988a,GrKo2004,N1980,N1986}.
The most important example of a Hilbert complex for our purposes of the
FEEC arises from the de~Rham complex of smooth differential forms on 
a domain or manifold.  

The main developments in FEEC to date have been for linear
(and now semi-linear) elliptic problems such as Poisson's equation
\begin{equation*}
-\Delta u = f.
\end{equation*}
Our goal here is to expand the scope of this analysis to include parabolic 
linear (and semi-linear) equations such as the heat equation,
\begin{equation*}
(\p_t-\Delta) u = f,
\end{equation*}
and hyperbolic equations such as the wave equation,
\begin{equation*}
(\p_{tt}-\Delta) u = f.
\end{equation*}
The exterior calculus framework treats $\Delta$ as $(d+\delta)^2$,
where $d$ is the exterior derivative operator and $\delta$ its adjoint.
The incorporation of the time derivative operation $\p_t$ into this framework,
however, has not been previously considered.
To remedy this, we develop the most natural extension of FEEC theory to
evolution problems: a generalization of the semi-discrete method
often called the `method of lines.'
This approach involves the discretization of the spatial part of
the differential operator, leaving the time variable continuous.
It can be viewed as introducing a time parameter into the discrete
(Hilbert) spaces that have been developed for elliptic problems.
The FEEC terminology provides a clear
and consistent notation for bounding errors in mixed methods
accumulated over a finite time interval.

We note that there is another approach to solving evolution problems
with finite elements, namely using a complete discretization of space-time.
This tactic allows for the dynamical
change of the underlying discrete approximation spaces in both space and time. 
Such an approach gives rise to space-time adaptivity, and is potentially
the most flexible and powerful approach to the numerical treatment of
parabolic and hyperbolic evolution problems.
This approach, which we will consider in a second article,
is most naturally formulated using geometric calculus, 
a well-studied mathematical structure for time-dependent problems.  
In the current article, we focus on extending FEEC to semi-discrete methods
using time-parameterized Banach (i.e.\ Bochner) space norm estimates 
for the method-of-lines approach.

Finally, we note that the work presented here was developed simultaneously and independently from a related project by Arnold and Chen~\cite{AC2012} for generalized Hodge-Laplacian style linear parabolic problems.
Our focus in this work is to extend the \textit{scalar} Hodge-Laplacian to both linear and semi-linear parabolic problems as well as linear hyperbolic problems, as this touches the existing literature on semi-discrete methods in the broadest fashion.
The pairing of these two results will lead to further insight in a variety of research directions.\\

{\em Outline and summary of contributions.}
The remainder of the paper is structured as follows.
In Section~\ref{sec:semidisc}, we review the classical semi-discrete 
mixed finite element method error estimates for parabolic problems (due to Thom\'ee~\cite{T2006} and others) and for hyperbolic problems (due to Geveci~\cite{G1988} and others).
In Section~\ref{sec:feec}, we give a very brief overview of the
Finite Element Exterior Calculus and recall some relevant results.
In Section~\ref{sec:results-par}, we combine the classical approach to semi-discrete
methods with modern FEEC theory
to establish some basic {\em a priori} error estimates for Galerkin mixed
finite element methods for parabolic problems.
The main result is Theorem~\ref{thm:boch-est-par}, which exploits the FEEC framework 
to obtain a classification of spatial finite element spaces that give optimal 
order convergence rates in Bochner norms.
In Section~\ref{sec:results-hyp}, we carry out a similar analysis for hyperbolic
problems, resulting in the error estimate given in Theorem~\ref{thm:boch-est-hyp}, a simultaneous sharpening of the result by Geveci for problems in two dimensional domains and a generalization to problems on $n$ dimensional domains.
Our results recover the estimates of Thom\'ee and Geveci 
for two-dimensional domains and the lowest-order mixed method as a special case, 
effectively extending their results to arbitrary spatial dimension and to
an entire family of mixed methods.
In Section~\ref{sec:semilin}, we employ the results of Holst and Stern~\cite{HS2010b} to
extend our parabolic estimates to a class of semi-linear evolution PDE, as encapsulated by Theorem~\ref{thm:semilin-est}.
Finally, in Section~\ref{sec:conc}, we draw conclusions and make 
remarks on future directions.
The appendices contain some technical details relating to the application of the semi-linear error estimates in the case $k=n$ and a review of some standard existence and uniqueness results for abstract parabolic and hyperbolic problems using the terminology of Bochner spaces.

\section{Semi-Discrete FEM Error Estimates for Evolution Problems}
\label{sec:semidisc}
We begin by reviewing semi-discrete finite element methods and their \textit{a priori} error estimates for parabolic and hyperbolic PDE systems.
We focus in each case on a relatively simple, well-studied system of interest to modeling communities, namely, the heat equation (parabolic) and the wave equation (hyperbolic). To avoid complicating the problem statements, we assume the given data $f$, $g$, $u_{0}$ and $u_{1}$ to be as regular as necessary for the problem context, unless otherwise specified.
The \textbf{heat equation} is: find $u(x,t)$ such that 
\begin{equation}
\label{eq:par-cnts}
\begin{tabular}{rlll}
$u_t - \Delta u = f$ & in $\Omega$, & for $t>0$ \\
$u = 0$ & on $\p\Omega$, & for $t>0$ & with $u(\cdot,0)=g$ in $\Omega$.
\end{tabular}
\end{equation}
We review the approach to Galerkin methods for this problem as presented in 
Thom\'ee~\cite{T2006} for domains $\Omega\subset\R^2$.  
His approach is based on work with Johnson~\cite{JT1981} and
builds upon prior analysis of elliptic projection~\cite{BSTW1977}.
A similar approach, restricted to $\Omega\subset \R^2$, was carried out by Garcia in~\cite{G1994}.
Let $\sigma=\nabla u$ and define the mixed, weak form problem: 
Given $f$ and $g$,
find $(u(t),\sigma(t)) \in L^{2}(\Omega)\times H(\Div; \Omega)$ such that 
\begin{equation}
\label{eq:par-mixedweak}
\begin{tabular}{rllll}
$(u_{t},\phi) - (\Div\, \sigma,\phi)$ & $= (f,\phi),$ & $\forall~ \phi\in L^2(\Omega),$ & $t>0$ \\
$ (u,\Div\,\omega) + (\sigma,\omega)$ & $= 0,$ & $\forall~\omega\in H(\Div;\Omega),$ & $t>0$, & $u(0)=g$.
\end{tabular}
\end{equation}
The semi-discrete problem is then to find $(u_h(t),\sigma_h(t))\in S_h\times H_h\subset L^2(\Omega) \times H(\Div;\Omega)$ such that
\begin{equation}
\label{eq:par-semidisc-first}
\begin{tabular}{rllll}
$(u_{h,t},\phi_h) - (\Div\, \sigma_h,\phi_h)$ & $= (f,\phi_h),$ & $\forall~ \phi_h\in S_h,$ & $t>0$ \\
$ (u_h,\Div\, \omega_h) + (\sigma_h,\omega_h)$ & $= 0,$ & $\forall~\omega_h\in H_h,$ & $t>0$, & $u_h(0)=g_h$.
\end{tabular}
\end{equation}
where $g_h$ is an approximation of $g$ in $S_h$.  
With bases for $S_h$ and $H_h$, the matrix form of the discrete problem is
\begin{align*}
AU_t - B\Sigma &= F, \\
B^TU + D\Sigma &= 0,\;\;\text{for $t>0$,\;\; $U(0)$ given},
\end{align*}
where $U$ and $\Sigma$ are vectors corresponding to $u_h$ and $\sigma_h$.
It is easily seen that the matrices $A$ and $D$ are positive definite.  Eliminating $\Sigma$, we have the system of ODEs
\[AU_t + BD^{-1}B^TU = F,\;\;\text{for $t>0$,\;\; $U(0)$ given},\]
which by standard results in ODE theory has a unique solution.

Thom\'ee uses discontinuous linear elements for $S_h$ and piecewise quadratic elements for $H_h$.  He defines the solution operator $T_h:L^2\raw S_h$ given by $T_h f = u_h$ for the corresponding elliptic problem and sets
\[g_h:=R_hg:=-T_h\Delta g.\]
For $g_h=R_hg$ and $t\geq 0$, Thom\'ee derives the estimates
\begin{align}
\vn{u_h(t)-u(t)}_{L^2} & \leq ch^2\left(\vn{u(t)}_{H^2} +\int_0^t\vn{u_t}_{H^2}ds \right), \label{eq:thomee-1-u}\\
\vn{\sigma_h(t)-\sigma(t)}_{L^2} & \leq ch^2\left(\vn{u(t)}_{H^3} +\left(\int_0^t\vn{u_t}_{H^2}^2ds\right)^{1/2} \right).\label{eq:thomee-1-s}
\end{align}
Note that these estimates are for a fixed time value $t$ and restricted to a particular choice of finite elements in 2D.

We now turn to the \textbf{wave equation}: find $u(x,t)$ such that
\begin{equation}
\label{eq:hyp-cnts}
\begin{tabular}{rlll}
$u_{tt} - \Delta u = f$ & in $\Omega$, & for $t>0$, \\
$u = 0$ & on $\p\Omega$, & for $t>0$ & with $u(\cdot,0)=u_0$ in $\Omega$, \\
 & & & and $u_t(\cdot,0) = u_1$ in $\Omega$
\end{tabular}
\end{equation}
There are two approaches to defining a mixed weak form of this problem.  The first is very similar to the parabolic case: 
find $(u(t),\sigma(t)) \in L^{2}(\Omega)\times H(\Div; \Omega)$ such that
\begin{equation}
\label{eq:hyp-mixedweak}
\begin{tabular}{rllll}
$(u_{tt},\phi) - (\Div~\sigma,\phi)$ & $= (f,\phi),$ & $\forall~ \phi\in L^2(\Omega),$ & $t>0$, \\[2mm]
$ (u,\Div~\omega) + (\sigma,\omega)$ & $= 0,$ & $\forall~\omega\in H(\Div; \Omega),$ & $t>0$,  \\[2mm]
$u(0)$ & $= u_0$, \\[2mm]
$u_t(0)$ & $= u_1$. \\[2mm]
\end{tabular}
\end{equation}
It is difficult to derive estimates for the numerical approximation of (\ref{eq:hyp-mixedweak}) akin to those found in the parabolic case due to the second derivatives appearing in the formulation.  
Some attempts at estimates along these lines for $\Omega\subset\R^2$ have been given by Baker~\cite{B1976} and Cowsar, Dupont and Wheeler~\cite{CDW1990,CDW1996}.

For the purpose of extending the FEEC framework, we find the \textbf{`velocity-stress' formulation} of the problem and the results of Geveci~\cite{G1988} to be more useful.
This formulation solves for $\ut:= u_t$ instead of $u$: Given $f$, $u_0$, and $u_1$, find $(\ut(t),\sigma(t)) \in L^{2}(\Omega)\times H(\Div; \Omega)$ such that
\begin{equation}
\label{eq:boch-hyper-vs}
\begin{tabular}{rllll}
$(\ut_t,\phi) - (\Div~\sigma,\phi)$ & $= (f,\phi),$ & $\forall~ \phi\in L^2(\Omega),$ & $t>0$, \\[2mm]
$(\ut,\Div~\omega)+(\sigma_t,\omega)$ & $= 0,$ & $\forall~\omega\in H(\Div; \Omega),$ & $t>0$, \\[2mm]
$\ut(0)$ & $= u_1$, \\[2mm]
$\sigma(0)$ & $= \nabla u_0$. \\[2mm]
\end{tabular}
\end{equation}
The semi-discrete problem is then to find $(\ut_h(t),\sigma_h(t)) \in S_{h}\times H_{h}$ such that 
\begin{equation}
\label{eq:hyp-semidisc-first}
\begin{tabular}{rllll}
$(\ut_{h,t},\phi_h) - (\Div~\sigma_h,\phi_h)$ & $= (f,\phi_h),$ & $\forall~ \phi_h\in S_h,$ & $t>0$, \\[2mm]
$(\ut_h,\Div~\omega_h) + (\sigma_{h,t},\omega_h)$ & $= 0,$ & $\forall~\omega_h\in H_h,$ & $t>0$, \\[2mm]
$\ut_h(0)$ & $=u_{1,h}$, \\[2mm]
$\sigma_h(0)$ & $=(\nabla u_0)_h$,
\end{tabular}
\end{equation}
where $u_{1,h}$ is an approximation of $u_1$ in $S_h$ and $(\nabla u_0)_h$ is an approximation of $\nabla u_0$.
Again, bases for $S_h$ and $H_h$ reduce the discrete problem to a matrix formulation:
\begin{align*}
AW_t - B\Sigma &= F, \\
B^TW + D\Sigma_t  & = 0,\;\;\text{for $t>0$,\;\; $W(0)$, $\Sigma(0)$ given},
\end{align*}
where $W$ and $\Sigma$ are vectors corresponding to $\ut_h$ and $\sigma_h$ and $A$ and $D$ are symmetric, positive definite matrices.
As Geveci~\cite[p. 248]{G1988} explains, this can be reduced to a single iterative system of the form
\[(D+k^2B^TA^{-1}B)\Sigma^{n+1}= G,\]
where $k$ denotes the time step in an implicit Euler time-differencing scheme.

To derive an error estimate for the velocity-stress discretization, Geveci states the need for projection operators from $\hdiv$ to $H_h$ and from $L^2$ to $S_h$ satisfying certain approximation properties.  
He explains that such operators exists for a variety of finite element spaces in $\R^2$, e.g.\ the Raviart-Thomas spaces~\cite{RT1977}, allowing the following result.
For $1\leq s\leq r$ with $r\geq 2$, 
\begin{equation}
\label{eq:geveci-est}
\begin{split}
\eltn{\ut_h(t)-\ut(t)}+ \eltn{\sigma_h(t)-\sigma(t)} \leq c\left(\eltn{u_1-u_{1,h}}+\eltn{\nabla u_0-(\nabla u_0)_h}\right) +\\
+ ch^s\left(\vn{u_1}_s+\vn{\nabla u_0}_s+\int_0^t(\vn{\ut_t(\tau)}_s+\vn{\sigma_t(\tau)}_s d\tau \right).
\end{split}
\end{equation}
Like estimates (\ref{eq:thomee-1-u}) and (\ref{eq:thomee-1-s}) for the parabolic problem, (\ref{eq:geveci-est}) says that the approximation error can be controlled in $L^2$ norm at any time $t$ by the $H^s$ norm of the initial conditions plus the accumulated norm of the variables up to time $t$.   It is these types of estimates that the FEEC framework can refine, simplify, and generalize to arbitrary spatial dimension $n$.

\section{The Finite Element Exterior Calculus}
\label{sec:feec}
The finite element exterior calculus (FEEC) provides an elegant mathematical
framework for deriving error estimates for a large class of elliptic PDE.
We now give a brief overview of the notation and certain  main results from FEEC 
which are relevant to this paper.  
We refer the reader to the seminal papers 
of Arnold, Falk, and Winther~\cite{AFW2006,AFW2010} for additional explanation.

Let $\Omega$ be a bounded $n$-manifold embedded in $\R^n$ and assume $\Omega$ has a piecewise smooth, Lipschitz boundary.  The space of $L^2$-bounded continuous differential $k$-forms on $\Omega$ is given by
\[L^2\Lambda^k(\Omega):=\left\{\sum_i a_i dx_i\in\Lambda^k(\Omega)\;:\;a_i\in L^2(\Omega)\quad\forall~ i\right\},\]
where $i$ ranges over all strictly increasing sequences of $k$ indices chosen from $\{1,\ldots,n\}$.  The exterior derivative operator $d_k:\Lambda^k(\Omega)\raw\Lambda^{k+1}(\Omega)$ acts on these spaces to form a Hilbert complex $(L^2\Lambda,d)$.  The associated domain complex is the sequence of spaces $H\Lambda^k:=\text{domain($d_k$)}\subset L^2\Lambda^k(\Omega)$, commonly called the $L^2$ deRham complex:
\[
\xymatrix{
0 \ar[r] & H\Lambda^0 \ar[r]^-{d_0} & H\Lambda^1 \ar[r]^-{d_1} & \cdots \ar[r]^-{d_{n-1}} & H\Lambda^n \ar[r] & 0.
}
\]
The norm on each space is the graph norm associated to $d$, i.e.
\[(u,v)_{H^k(\Omega)}:=(u,v)_{\Lambda^k(\Omega)}+(d_ku,d_kv)_{\Lambda^{k+1}(\Omega)}.\]
We note that in any dimension $n$, the beginning and end of the $L^2$ deRham complex can be understood in terms of traditional Sobolev spaces and differential operators:
\[
\xymatrix{
0 \ar[r] & H^1(\Omega) \ar[r]^-{\text{grad}} &  \cdots \ar[r]^-{\Div} & L^2(\Omega) \ar[r] & 0
}
\]
A major conclusion of FEEC is that stable finite element methods for elliptic PDE must seek solutions in finite dimensional subspaces $\Lambda^k_h\subset H\Lambda^k$ that satisfy certain key approximation properties.  
First, the subspaces should form a subcomplex of the $L^2$ deRham complex, meaning $d\Lambda^k_h\subset\Lambda^{k+1}_h$.  
Second, $\Lambda^k_h$ should have sufficient approximation that upper bounds on $\inf_{v\in \Lambda^k_h}\vn{u-v}_{H\Lambda^k}$ can be ensured for some or all $u\in H\Lambda^k$.  
Third, there must exist bounded cochain projections $\pi_h^k:H\Lambda^k\raw\Lambda^k_h$ which are invariant on $\Lambda^k_h$, commute with the exterior derivative operators, and provide a bound $\vn{\pi^k_hv}_{H\Lambda^k}\leq c\vn{v}_{H\Lambda^k}$ for all $v\in{H\Lambda^k}$.

In the context of the deRham complex, all these properties are shown to be provided for by two canonical classes of piecewise degree $r$ polynomials associated to a simplicial mesh $\Tau$ of $\Omega$.  
Let $\poly_r$ denote polynomials in $n$ variables of degree at most $r$ and $\mathcal H_r\subset\poly_r$ the subspace of homogeneous polynomials.   
The first class, denoted $\poly_r\Lambda^k(\Tau)$, consists of all $k$-forms with coefficients belonging to $\poly_r$ on each $n$-simplex of $\Tau$.  
The second class, denoted $\poly_r^-\Lambda^k(\Tau)$, interleaves with the first class, i.e.
\[\poly_{r-1}\Lambda^k(\Tau)\subsetneq \poly_r^-\Lambda^k(\Tau) \subsetneq \poly_r\Lambda^k(\Tau).\]
To define $\poly_r^-\Lambda^k(\Tau)$, first define $X$ be the vector field on $\R^n$ such that $X(x)$ is the vector based at $x\in\R^n$ that points opposite to the origin with length $|x|$.
Define $\poly_r^-\Lambda^k(\Tau):=\poly_r\Lambda^k\oplus\kappa\mathcal H_{r-1}\Lambda^{k+1}$, a direct sum, where $\kappa$ is defined by contraction with $X$.  
The map $\kappa$ is called the \textit{Koszul differential} and gives rise to the \textit{Koszul complex}.
This is elaborated upon in detail in the work of Arnold, Falk and Winther~e.g.~\cite[p. 328]{AFW2010}.

For $n=3$, we have the following correspondences between the FEEC notation of finite element spaces and traditional finite element spaces.
\begin{center}
\begin{tabular}{rcl}
$\poly_{r+1}\Lambda^2(\Tau)$ &=& \nedelec 2nd-kind $\hdiv$ elements of degree $\leq r+1$ (see~\cite{N1986})\\[2mm]
$\poly_{r+1}^-\Lambda^2(\Tau)$ &=& \nedelec 1st-kind $\hdiv$ elements of order $r+1$ (see~\cite{N1980})\\[2mm]
$\poly_{r+1}^-\Lambda^3(\Tau)$ &=& $\poly_{r}\Lambda^3(\Tau)$\quad=\quad discontinuous elements of degree $\leq r$ 
\end{tabular}
\end{center}
Hence, in the case of the deRham complex, FEEC recovers well-known finite element spaces while at the same time describing their generalization to arbitrary spatial dimensions.

The last piece of FEEC used in this work is the existence of smoothed projection operators
\begin{equation}
\label{eq:coch-proj-def}
\pi_h^k:L^2\Lambda^k\raw L^2\Lambda^k_h\quad\text{where~$\Lambda^k_h\in\{\poly_r\Lambda^k(\Tau),~\poly_r^-\Lambda^k(\Tau)\}$}.
\end{equation}
These operators are shown, by virtue of their construction, to be uniformly bounded (in $L^2\Lambda^k$) with respect to $h$. 
An explicit construction of these operators can be found in the papers of Arnold, Falk, and Winther~\cite{AFW2006,AFW2010}; the following theorem asserts some of their key properties.
\begin{theorem}[\cite{AFW2010} Theorem 5.9]\hfill
\label{thm:afw-5pt9}
\renewcommand{\labelenumi}{(\roman{enumi}.)}
\begin{enumerate}
\item Let $\Lambda^k_h$ be one of the spaces $\poly^-_{r+1}\Lambda^k(\Tau)$ or, if $r\geq 1$, $\poly_r\Lambda^k(\Tau)$.  Then $\pi^k_h$ is a projection onto $\Lambda^k_h$ and satisfies
\[\vn{\omega-\pi^k_h\omega}_{L^2\Lambda^k(\Omega)}\leq ch^s\vn{\omega}_{H^s\Lambda^k(\Omega)},\quad \omega\in H^s\Lambda^k(\Omega),\]
for $0\leq s\leq r+1$.  Moreover, for all $\omega\in L^2\Lambda^k(\Omega)$, $\pi^k_h\omega\raw\omega$ in $L^2$ as $h\raw 0$.
\item Let $\Lambda^k_h$ be one of the spaces $\poly_r\Lambda^k(\Tau)$ or $\poly^-_r\Lambda^k(\Tau)$ with $r\geq 1$.  Then 
\[\vn{d(\omega-\pi^k_h\omega)}_{L^2\Lambda^k(\Omega)}\leq ch^s\vn{d\omega}_{H^s\Lambda^k(\Omega)},\quad \omega\in H^s\Lambda^k(\Omega),\]
for $0\leq s\leq r$.
\item Let $\Lambda^{k-1}_h\in\left\{\poly_{r+1}\Lambda^{k-1}(\Tau), \poly_{r+1}^-\Lambda^{k-1}(\Tau)\right\}$ and $\Lambda^k_h= \poly_{r+1}^-\Lambda^k(\Tau)$ or, if $r>0$, $\poly_r\Lambda^k(\Tau)$.  Then $d\pi^{k-1}_h=\pi^k_hd$.
\end{enumerate}
\end{theorem}

To handle the setting of parabolic and hyperbolic PDE,
we now merge the notation of FEEC with the abstract framework of
parametrized Banach spaces, also known as Bochner spaces.
We follow prior approaches along these lines, especially~\cite[page 66]{T1997} and \cite{RR2004}.
Let $X$ be a Banach space and $\mr I:=(0,T)$ an interval with closure $I := \ol I =[0,T]$.  
Define the space
\[C(\mr I,X) := \{u:\mr I\raw X\;\;|\;\;\text{$u$ bounded and continuous}\},\]
and equip it with the norm
\[\vn{u}_{C(\mr I,X)} := \sup_{t\in \mr I}\vn{u(t)}_X.\]
The \textbf{Bochner space} $L^P(\mr I,X)$ is then defined to be the completion of $C(\mr I,X)$ with respect to the norm
\[\vn{u}_{L^p(\mr I,X)} := \left(\int_{\mr I}\vn{u(t)}^p_X dt\right)^{1/p}.\]
The space $H^1(\mr I,X)$ has an analogous norm
\[\vn{u}_{H^1(\mr I,X)} := \left(\int_{\mr I}\vn{u(t)}^2_X + \vn{\frac d{dt}u(t)}^2_X dt\right)^{1/2}.\]
We will commonly use $X=L^2\Lambda^k$ or $X=H^s\Lambda^k$ where it is understood that the forms are defined over a spatial domain $\Omega$. In the case of spatial-only norms, we will use the notation $\vn{\cdot}_{L^2}$ or just $\vn{\cdot}$ for the $L^{2}\Lambda^{k}$ norm, and $\|\cdot\|_{H^{s}}$ for the $H^{s}\Lambda^k$ norm, where $k$ will be clear from context.

\section{\emph{A Priori} Error Estimates for Parabolic Problems}
\label{sec:results-par}
We extend Thom\'ee's error estimates from Section~\ref{sec:semidisc} to the broader class of elements and arbitrary spatial dimension allowed by FEEC with Bochner space norms.
For simplicity, let $\Omega\subset\R^n$ be a contractible domain.
Define the \textbf{mixed weak parabolic} problem: 
Given $f,g$, find~$(u,\sigma):I\raw H\Lambda^n\times H\Lambda^{n-1}$ such that
\begin{equation}
\label{eq:boch-mixedweak}
\begin{tabular}{rllll}
$(u_{t},\phi) - (\Div~\sigma,\phi)$ & $= (f,\phi),$ & $\forall~ \phi\in H\Lambda^n,$ & $t\in I$, \\[2mm]
$ (u,\Div~\omega)+ (\sigma,\omega)$ & $= 0,$ & $\forall~\omega\in H\Lambda^{n-1},$ & $t\in I$,  \\[2mm]
$u(0)$ & $= g$.
\end{tabular}
\end{equation}
The \textbf{semi-discrete parabolic problem} is thus: Find $(u_h,\sigma_h):I\raw \Lambda^n_h\times \Lambda^{n-1}_h$ such that
\begin{equation}
\label{eq:par-semidisc}
\begin{tabular}{rllll}
$(u_{h,t},\phi_h) - (\Div~\sigma_h,\phi_h)$ & $= (f,\phi_h),$ & $\forall~ \phi_h\in \Lambda^n_h,$ & $t\in I$, \\[2mm]
$ (u_h,\Div~\omega_h)+ (\sigma_h,\omega_h)$ & $= 0,$ & $\forall~\omega_h\in \Lambda^{n-1}_h,$ & $t\in I$,\\[2mm]
$u_h(0)$ & $=g_h$.
\end{tabular}
\end{equation}
Define $g_h$ to be the solution to the elliptic problem with load data $-\Delta g$, i.e.
\begin{equation}
\label{eq:gh-def}
\begin{tabular}{rllll}
$(\Div~\hat\sigma_h,\phi_h)-(\Delta g,\phi_h)$ & $=0$ & $\forall~ \phi_h\in \Lambda^n_h,$\\[2mm]
$(\hat\sigma_h,\omega_h) + (g_h,\Div~\omega_h)$ & $= 0,$ & $\forall~\omega_h\in \Lambda^{n-1}_h.$ 
\end{tabular}
\end{equation}
It is shown in~\cite{T2006} that a unique solution to (\ref{eq:par-semidisc}) exists, based on the positive-definiteness of the solution operator $T_h:L^2\raw \Lambda^n_h$ for the elliptic problem.
A more basic argument for this result can also be made by appealing to the existence of an adjoint to the discrete divergence operator.

Elliptic projection, an idea dating back to Wheeler~\cite{W1973}, can be carried out for any fixed time value as we now discuss.  For any $t_0\in I$, define the \textbf{time-ignorant discrete elliptic problem}: Find $(\tilde u_h,\tilde\sigma_h)\in \Lambda^n_h\times\Lambda^{n-1}_h$ such that
\begin{equation}
\label{eq:par-to-ellip}
\begin{tabular}{rll}
$(\Div~\tilde\sigma_h,\phi_h) + (-\Delta u(t_0),\phi_h)$ & $= 0,$ & $\forall~ \phi_h\in \Lambda^n_h,$  \\[2mm] 
$(\tilde\sigma_h,\omega_h) + (\tilde u_h,\Div~\omega_h)$ & $= 0,$ & $\forall~\omega_h\in \Lambda^{n-1}_h,$  \\[2mm] 
$\tilde u_h(0)$ & $= g_h$.
\end{tabular}
\end{equation}
Note that the $u$ appearing in the first equation of (\ref{eq:par-to-ellip}) is the solution to the continuous problem (\ref{eq:boch-mixedweak}).  Thus, we can view $\tilde\sigma_h$ and $\tilde u_h$ as functions of $t$ with the understanding that they are defined for each $t$ value by (\ref{eq:par-to-ellip}) alone; no continuity with respect to $t$ is required, hence the moniker `time-ignorant.'   

For ease of notation, and in keeping with Thom\'ee, define the error functions
\begin{align*}
\rho(t) & := \tilde u_h(t) - u(t), \\
\theta(t) & := u_h(t) -\tilde u_h(t), \\
\veps(t) & := \sigma_h(t) -\tilde\sigma_h(t).
\end{align*}
We now prove a lemma which will aid in our subsequent analysis.
The result appears as part of the proof of Thom\'ee~\cite[Theorem 17.2]{T2006} but we expand it here for clarity.
\begin{lemma}[Thom\'ee~\cite{T2006}]
\label{lem:err-eqns-lin}
The error functions satisfy the semi-discrete formulation:
\begin{equation}
\label{eq:par-error}
\begin{tabular}{rllll}
$(\theta_{t},\phi_h) - (\text{\em div}~\veps,\phi_h)$ & $= -(\rho_t,\phi_h),$ & $\forall~ \phi_h\in \Lambda^n_h,$ & $t\in I$, \\[2mm]
$ (\theta,\text{\em div}~\omega_h)+ (\veps,\omega_h)$ & $= 0,$ & $\forall~\omega_h\in \Lambda^{n-1}_h,$ & $t\in I$.
\end{tabular}
\end{equation}
\end{lemma}
\begin{proof}
The second equation is immediate from the second equations in (\ref{eq:par-semidisc}) and (\ref{eq:par-to-ellip}).  The first equation can be written out as
\begin{align*}
(u_{h,t},\phi_h) - (\Div~\sigma_h,\phi_h) + (\Div~\tilde\sigma_h,\phi_h) - (\tilde u_{h,t},\phi_h) & =  (u_t,\phi_h) - (\tilde u_{h,t},\phi_h)
\end{align*}
which is reduced as follows:
\begin{align*}
(u_{h,t},\phi_h) - (\Div~\sigma_h,\phi_h) + (\Div~\tilde\sigma_h,\phi_h)  & =  (u_t,\phi_h) & \text{cancel like terms} \\
(u_{h,t},\phi_h) - (\Div~\sigma_h,\phi_h) & = -(\Delta u,\phi_h)+  (u_t,\phi_h) & \text{by (\ref{eq:par-to-ellip})} \\
(f,\phi_h)  & = -(\Delta u,\phi_h)+  (u_t,\phi_h) & \text{ by (\ref{eq:par-semidisc})}
\end{align*}
This says that the continuous problem $u_t-\Delta u = f$ should hold in a weak sense when tested against any of the functions in $\Lambda^n_h$.  This is guaranteed to be true since we chose $\Lambda^n_h\subset\Lambda^n=L^2$.  Thus, the error equations hold as stated.
\end{proof}

The following theorem says that if $\Lambda^n_h$ and $\Lambda^{n-1}_h$ are chosen according to the FEEC framework, then error estimates akin to (\ref{eq:thomee-1-u}) and (\ref{eq:thomee-1-s}) can be obtained.
Note that in the semidiscrete setting, $(\Delta u)_t(t) = \p_t\Delta u(t)$ since the time and spatial derivatives commute, allowing the simplified notation $\Delta u_t(t)$ used here.

\begin{theorem}
\label{thm:boch-est-par}
Fix $\Omega\subset\R^n$ and fix $I:=[0,T]$.  
Suppose $(u,\sigma)$ is the solution to (\ref{eq:boch-mixedweak}) such that the regularity estimate
\begin{equation}
\label{eq:reg-est}
\vn{u(t)}_{H^{s+2}}+\vn{\sigma(t)}_{H^{s+1}}\leq c\vn{\Delta u(t)}_{H^{s}}
\end{equation}
holds for $0\leq s\leq s_{\max}$ and $t\in I$. 

Choose finite element spaces 
\[\Lambda^{n-1}_h = \left\{\begin{array}{c} \poly_{r+1}\Lambda^{n-1}(\Tau) \\ \text{or}\\ \poly_{r+1}^-\Lambda^{n-1}(\Tau)\end{array}\right\},\quad \Lambda^n_h = \poly_{r+1}^-\Lambda^n(\Tau) \;\;\left(= \poly_{r}\Lambda^n(\Tau)\right)
\]
Then for $0\leq s\leq s_{\max}$, $g_h$ defined by (\ref{eq:gh-def}), and $(u_h,\sigma_h)$ the solution to (\ref{eq:par-semidisc}), the following error estimates hold:
\begin{align}
\vn{u_h-u}_{L^2(I,L^2\Lambda^n)} &\leq \left\{\begin{array}{ll} \ds ch\left(\vn{\Delta u}_{L^2(I,L^2)}+\sqrt T\vn{\Delta u_t}_{L^1(I,L^2)}\right) & \text{if $r=0$} \\[4mm] \ds ch^{2+s}\left(\vn{\Delta u}_{L^2(I,H^s)}+ \sqrt T\vn{\Delta u_t}_{L^1(I,H^s)}\right) & \text{for $r>0$,} \\
& \text{~~if $s\leq r-1$}\end{array}\right.\label{eq:new-type1-u}
\end{align}
\begin{align}
\vn{\sigma_h-\sigma}_{L^2(I,L^2\Lambda^{n-1})} & \leq \left\{\begin{array}{l} \ds ch\left(\vn{\Delta u}_{L^2(I,H^s)}+\sqrt T\vn{\Delta u_t}_{L^2(I,L^2)}\right) \\[3mm] \qquad \text{if~$r=0$, $s=0$, $\Lambda^{n-1}_h=\poly_1^-\Lambda^{n-1}(\Tau)$} \\[3mm]
\ds c\left(h^{1+s}\vn{\Delta u}_{L^2(I,H^s)}+h\sqrt T\vn{\Delta u_t}_{L^2(I,L^2)}\right) \\[3mm] \qquad\text{if~$r=0$, $s\leq 1$, $\Lambda^{n-1}_h=\poly_1\Lambda^{n-1}(\Tau)$} \\[3mm]
\ds c\left(h^{1+s}\vn{\Delta u}_{L^2(I,H^s)}+h^{(3/2)+s}\sqrt T\vn{\Delta u_t}_{L^2(I,H^s)}\right) \\[3mm] \qquad\text{for $r>0$, if $s\leq r-1$}
\end{array}\right. \label{eq:new-type1-s}
\end{align}
\begin{align}
\vn{\textnormal{div}(\sigma_h-\sigma)}_{L^2(I,L^2\Lambda^n)} & \leq \left\{\begin{array}{l} \ds c\left(h^s\vn{\Delta u}_{L^2(I,H^s)}+h\vn{\Delta u_t}_{L^2(I,L^2)}\right) \\[3mm] \qquad \text{if~$r=0$, $s\leq 1$}\\[3mm]
\ds c\left(h^s\vn{\Delta u}_{L^2(I,H^s)}+h^{2+s}\vn{\Delta u_t}_{L^2(I,H^s)}\right) \\[3mm] \qquad\text{for $r>0$, if $s\leq r-1$}
\end{array}\right.\label{eq:new-type1-ds}
\end{align}
\end{theorem}

\begin{remark}
\em
Previous literature on semi-discrete methods usually leaves regularity assumptions implied by the error estimates.
For instance, if $\vn{u(t)}_{H^3}$ appears on the right side, it is implicitly assumed that $u(t)\in H^3$ for all $t\in I$.
We have stated the specific regularity assumption~(\ref{eq:reg-est}) to make clear what regularity must be assumed and to follow the presentation from Arnold, Falk, and Winther~\cite[p. 342]{AFW2010}.
The careful reader will notice that the left side of~(\ref{eq:reg-est})  does not include a $du(t)$ term, since $u(t)\in \Lambda^n$ implies $du(t)=0$, nor a $\vn{d\sigma(t)}_{H^{s}}$ term, since this can be absorbed into the $\vn{\sigma(t)}_{H^{s+1}}$ term.
An additional, more subtle difference is that $f$ on the right side of~(\ref{eq:reg-est}) has been replaced by $\Delta u$.
While these two are equivalent in the elliptic case, $\Delta u(t)$ evolves based on the initial data $g$ while $f(t)$ is prescribed, meaning they are in general different in the parabolic setting.
\end{remark}

\begin{proof}
Observe that (\ref{eq:par-to-ellip}) is exactly the $k=n$ case of the Hodge-Laplacian problem analyzed by Arnold, Falk and Winther~\cite{AFW2010} and the hypotheses here match their hypotheses.  We can thus use a triangle inequality argument for each estimate, e.g.
\begin{equation}
\label{eq:u-err-tri-ineq}
\vn{u(t)-u_h(t)} \leq \vn{u(t)-\tilde u_h(t)} + \vn{\tilde u_h(t)-u_h(t)}= \vn{\rho(t)} + \vn{\theta(t)}
\end{equation}
The first term will be bounded using the estimates from~\cite{AFW2010} and the second by the techniques from Thom\'ee~\cite{T2006}.   
The FEEC estimate~\cite[p. 342]{AFW2010} gives immediately 
\begin{equation}
\label{eq:rho-est}
\|\rho(t) \|\leq\left\{\begin{array}{ll} ch \|\Delta u(t)\| & \text{if $r=0$} \\[2mm] ch^{2+s}\vn{\Delta u(t)}_{H^s} & \text{if $s\leq r-1$, for $r>0$}\end{array}\right.
\end{equation}
Bounding $\eltn{\theta(t)}$ is more subtle.  Set $\phi_h := \theta$ and $\omega_h := \veps$ in (\ref{eq:par-error}).  Adding the equations yields
\begin{equation}
\label{eq:norms-to-ip}
\frac 12 \frac{d}{dt}\vn{\theta}^2 + \vn{\veps}^2 = -(\rho_t,\theta),\;\; t\in I
\end{equation}
We use a technique from Thom\'ee~\cite[p. 8]{T2006} to derive an estimate for $\vn{\theta(t)}$.  Since $\vn{\theta}$ may not be differentiable when $\theta=0$, introduce a constant $\delta>0$ and observe that
\[(\vn{\theta}^2+\delta^2)^{1/2}\frac d{dt}(\vn{\theta}^2+\delta^2)^{1/2}=\frac 12\frac{d}{dt}(\vn{\theta}^2+\delta^2)=\frac 12\frac{d}{dt}\vn{\theta}^2\leq\vn{\rho_t}\vn{\theta},\]
the last step following by (\ref{eq:norms-to-ip}) and Cauchy-Schwarz.  Since $\vn{\theta}\leq (\vn{\theta}^2+\delta^2)^{1/2}$, we have that
\[\frac{d}{dt}(\vn{\theta}^2+\delta^2)^{1/2}\leq\vn{\rho_t}.\]
Note that $\theta(0)=u_h(0)-\tilde u_h(0)=g_h-g_h=0$.  Thus
\[
\vn{\theta(t)}  = \lim_{\delta\raw 0}\int_0^t\frac{d}{dt}(\vn{\theta}^2+\delta^2)^{1/2}\leq\int_0^t\vn{\rho_t}.
\]
Using the bounds on $\vn{\rho(t)}$ from (\ref{eq:rho-est}), we get
\begin{equation}
\label{eq:theta-est}
\vn{\theta(t)}\leq \left\{\begin{array}{ll} \ds ch\int_0^t \|\Delta u_t(\ell)\| d\ell & \text{if $r=0$,} \\[4mm] \ds ch^{2+s}\int_0^t\vn{\Delta u_t(\ell)}_{H^s}d\ell & \text{for $r>0$, if $s\leq r-1$.}\end{array}\right.
\end{equation}
We can now assemble estimate (\ref{eq:new-type1-u}) by collecting our results.  We show the technique of the case $r=0$ as the other case employs identical analysis.
\begin{align*}
\vn{u_h-u}_{L^2(I,L^2\Lambda^n)} & = \left(\int_0^T\vn{u_h(t)-u(t)}^2dt\right)^{1/2} \\
 & \leq \left(\int_0^T\left(\vn{\rho(t)}+\vn{\theta(t)}\right)^2dt\right)^{1/2} \\
 & \leq ch\left(\int_0^T\left(\vn{\Delta u(t)}+\int_0^t \|\Delta u_t(\ell)\| d\ell\right)^2dt\right)^{1/2} \\
 & \leq ch\left(\int_0^T 2\left(\vn{\Delta u(t)}^2+\left(\int_0^t \|\Delta u_t(\ell) \| d\ell\right)^2\right)dt\right)^{1/2}.
\end{align*}
Roll the 2 into the constant $c$ and observe that the inner integral is maximal when $t=T$.  Thus,
\begin{align*}
\vn{u_h-u}_{L^2(I,L^2\Lambda^n)} 
 & \leq ch\left(\int_0^T \vn{\Delta u(t)}^2+\vn{\Delta u_t}_{L^1(I,L^2)}^2 dt\right)^{1/2} \\
 &  =  ch\left(\vn{\Delta u}_{L^2(I,L^2)}^2 + T \vn{\Delta u_t}_{L^1(I,L^2)}^2 \right)^{1/2} \\
 & \leq ch\left(\vn{\Delta u}_{L^2(I,L^2)}^2 + T \vn{\Delta u_t}_{L^1(I,L^2)}^2 + \right.\\
 & \qquad\qquad \left. 2\vn{\Delta u}_{L^2(I,L^2)}\sqrt T \vn{\Delta u_t}_{L^1(I,L^2)}\right)^{1/2} \\
 & = ch\left(\vn{\Delta u}_{L^2(I,L^2)} + \sqrt T\vn{\Delta u_t}_{L^1(I,L^2)}\right).
\end{align*}

We now turn to (\ref{eq:new-type1-s}), i.e.\ an error bound for the approximation of $\sigma$.  We use the same technique of bounding $\vn{\sigma(t)-\tilde\sigma_h(t)}$ by the corresponding FEEC estimate and $\vn{\tilde\sigma_h(t)-\sigma_h(t)}$ ($=\vn{\veps(t)}$) by a modification of (\ref{eq:par-error}).  First, observe that the FEEC estimate~\cite[p. 342]{AFW2010} gives
\begin{equation}
\label{eq:AFWsigma-est}
\vn{\sigma(t)-\tilde\sigma_h(t)}\leq ch^{1+s}\vn{\Delta u(t)}_{H^s},~\text{if}~\left\{\begin{array}{ll} s\leq r+1, & \Lambda^{n-1}_h=\poly_{r+1}\Lambda^{n-1}(\Tau) \\[2mm] 
s\leq r, & \Lambda^{n-1}_h=\poly_{r+1}^-\Lambda^{n-1}(\Tau) \end{array}\right..
\end{equation}
To bound $\vn{\eps(t)}$, differentiate the second equation of (\ref{eq:par-error}) with respect to $t$ and set $\phi_h := 2\theta_t$, $\omega_h := 2\veps$, yielding
\begin{equation*}
\begin{tabular}{rllll}
$(\theta_{t},2\theta_t) - (\Div~\veps,2\theta_t)$ & $= -(\rho_t,2\theta_t),$\\[2mm]
$(\veps_t,2\veps) + 2(\theta_t,\Div~\veps)$ & $= 0.$ 
\end{tabular}
\end{equation*}
Adding the equations and converting to norms, we have the bound
\[\frac {d}{dt} \vn{\veps}^2+ 2\vn{\theta_t}^2 = -2(\rho_t,\theta_t) \leq \vn{\rho_t}^2 +\vn{\theta_t}^2,\]
by Cauchy-Schwarz and the AM-GM inequality.  Note that since $\theta(0)=0$, we have $\veps(0)=0$ by the second equation of (\ref{eq:par-error}).  Thus
\begin{equation}
\label{eq:eps-lin}
\vn{\veps(t)}^2 = \int_0^t\frac {d}{ds}\vn{\veps(s)}^2ds \leq \int_0^t \frac {d}{ds}\vn{\veps(s)}^2 +\vn{\theta_t}^2 ds \leq \int_0^t \vn{\rho_t}^2 ds.
\end{equation}
As before, we use (\ref{eq:rho-est}) to derive
\begin{equation}
\label{eq:veps-est}
\vn{\veps(t)} \leq 
\left\{\begin{array}{ll} \ds ch\left(\int_0^t\|\Delta u_t(\ell)\| d\ell\right)^{1/2} & \text{if $r=0$,} \\[4mm] \ds ch^{2+s}\left(\int_0^t\vn{\Delta u_t(\ell)}_{H^s}^2 d\ell\right)^{1/2} & \text{for $r>0$, if $s\leq r-1$.}\end{array}\right.
\end{equation}
We assemble estimate (\ref{eq:new-type1-s}) by collecting these results in the same fashion we did for (\ref{eq:new-type1-u}).  In the $r=0$ cases, we have
\begin{align*}
\vn{\sigma_h-\sigma}_{L^2(I,L^2\Lambda^{n-1})} & = \left(\int_0^T\vn{\sigma_h(t)-\sigma(t)}^2dt\right)^{1/2} \\
 & \leq \left(\int_0^T\left(\vn{\sigma(t)-\tilde\sigma_h(t)} +\vn{\veps(t)}\right)^2dt\right)^{1/2} \\
 & \leq c\left(\int_0^T\left(h^{1+s}\vn{\Delta u(t)}_{H^s}+h\left(\int_0^t\|\Delta u_t(\ell)\| d\ell\right)^{1/2} \right)^2dt\right)^{1/2} \\
 & \leq ch\left(\int_0^T 2\left(h^{2s}\vn{\Delta u(t)}_{H^s}^2+\int_0^t \|\Delta u_t(\ell)\| d\ell\right)dt\right)^{1/2}.
\end{align*}
Rolling the 2 into the constant $c$ and again noting that the inner integral is maximal when $t=T$, we recover the first two estimates of (\ref{eq:new-type1-s}):
\begin{align*}
\vn{\sigma_h-\sigma}_{L^2(I,L^2\Lambda^{n-1})}
 & \leq ch\left(\int_0^T h^{2s}\vn{\Delta u(t)}_{H^s}^2 +\vn{\Delta u_t}_{L^2(I,L^2)}^2 dt\right)^{1/2} \\
 &    = ch\left(h^{2s}\vn{\Delta u}_{L^2(I,H^s)}^2 + T\vn{\Delta u_t}_{L^2(I,L^2)}^2 dt\right)^{1/2} \\
 & \leq ch\left(h^s\vn{\Delta u}_{L^2(I,H^s)}+\sqrt T\vn{\Delta u_t}_{L^2(I,L^2)}\right).
\end{align*}
When $r>0$, (\ref{eq:AFWsigma-est}) requires $s\leq r$ or $s\leq r+1$ to obtain optimal convergence rates on the first term of the right side while (\ref{eq:veps-est}) requires $s\leq r-1$ to obtain optimal rates on the second term.  Thus the hypothesis $s\leq r-1$ implies both (\ref{eq:AFWsigma-est}) and (\ref{eq:veps-est}).
The last estimate of (\ref{eq:new-type1-s}) then follows by identical analysis to the first two cases.

Finally, we turn to estimate (\ref{eq:new-type1-ds}) and follow the same technique.  Since $\Div$ is the exterior derivative $d:\Lambda^{n-1}\rightarrow\Lambda^n$, we have from FEEC~\cite[p. 342]{AFW2010} that
\begin{equation}
\vn{\Div(\sigma(t)-\tilde\sigma_h(t))}\leq ch^{s}\vn{\Delta u(t)}_{H^s},~\text{if~$s\leq r+1$}.
\end{equation}
Taking the $L^2$ norm with respect to time over $[0,T]$ yields
\begin{equation}
\label{eq:AFWdsigma-est}
\vn{\Div(\sigma-\tilde\sigma_hŒ)}_{L^2(I,L^2)} \leq ch^{s}\vn{\Delta u}_{L^2(I,H^s)},~\text{if~$s\leq r+1$}.
\end{equation}
To bound $\vn{\Div~\veps}$ set $w_h:=\veps$ in~(\ref{eq:par-error}) and take the derivative with respect to $t$.  This yields
\[\frac{d}{dt}\vn{\veps}^2+(\theta_t,\Div~\veps)=0.\]
Note that $\Div~\veps\in \Lambda^n_h$ since the discrete spaces are chosen to satisfy the relationship $\Div\Lambda^{n-1}_h\subset\Lambda^n_h$.  Thus, we can set $\phi_h :=\Div~\veps$ in in~(\ref{eq:par-error}) and substitute to get
\begin{equation}
\label{eq:for-div-eps-est}
\frac{d}{dt}\vn{\veps}^2+\vn{\Div~\veps}^2 = (\rho_t,\Div~\veps).
\end{equation}
Integrating both sides over $[0,T]$ and applying Young's inequality, we have
\begin{equation}
\vn{\veps(T)}^2- \vn{\veps(0)}^2 + \int_0^T\vn{\Div~\veps}^2 dt  \leq \frac 12\int_0^T\vn{\rho_t}^2 + \vn{\Div~\veps}^2 dt
\end{equation}
Note that $\veps(0)=\sigma_h(0)-\tilde\sigma_h(0)=0$ since the elliptic projection is the identity at time $t=0$.
We thus have
\begin{equation}
\int_0^T\vn{\Div~\veps}^2 dt  \leq \int_0^T\vn{\rho_t}^2 dt.
\end{equation}
and by taking square roots,
\begin{equation}
\vn{\Div~\veps}_{L^2(I,L^2)}  \leq  \vn{\rho_t}_{L^2(I,L^2)} dt.
\end{equation}
Again, we use (\ref{eq:rho-est}) to get
\begin{equation}
\label{eq:div-veps-est}
\vn{\Div~\veps}_{L^2(I,L^2)}\leq \left\{\begin{array}{ll} \ds ch\vn{\Delta u_t}_{L^2(I,L^2)} & \text{if $r=0$} \\[2mm]\ds  ch^{2+s}\vn{\Delta u_t}_{L^2(I,H^s)} & \text{if $s\leq r-1$, for $r>0$}\end{array}\right..
\end{equation}
The estimate (\ref{eq:new-type1-ds}) follows by combining this with (\ref{eq:AFWdsigma-est}).
\end{proof}

\section{\emph{A Priori} Error Estimates for Hyperbolic Problems}
\label{sec:results-hyp}
We now analyze hyperbolic problems again using the Bochner space norms introduced at the end of Section~\ref{sec:feec}.
Again, let $\Omega\subset\R^n$ be a contractible domain.
Define the \textbf{velocity-stress mixed weak formulation}: 
Given the data functions $f, u_0, u_1$,
find $(\ut,\sigma): I\raw  H\Lambda^n \times H\Lambda^{n-1}$ such that
\begin{equation}
\label{eq:boch-hyper-s}
\begin{tabular}{rllll}
$(\ut_t,\phi) - (\Div~\sigma,\phi)$ & $= (f,\phi),$ & $\forall~ \phi\in H\Lambda^n,$ & $t\in I$, \\[2mm]
$(\ut,\Div~\omega) + (\sigma_t,\omega)$ & $= 0,$ & $\forall~\omega\in H\Lambda^{n-1},$ & $t\in I$, \\[2mm]
$\ut(0)$ & $= u_1$, \\[2mm]
$\sigma(0)$ & $= \nabla u_0$, \\[2mm]
\end{tabular}
\end{equation}
where $\ut = u_t$ as in (\ref{eq:boch-hyper-vs}).  
The \textbf{semi-discrete hyperbolic problem} is thus: Find $(\ut_h,\sigma_h):I\raw \Lambda^n_h\times \Lambda^{n-1}_h$ such that
\begin{equation}
\label{eq:hyp-semidisc}
\begin{tabular}{rllll}
$(\ut_{h,t},\phi_h) - (\Div~\sigma_h,\phi_h)$ & $= (f,\phi_h),$ & $\forall~ \phi_h\in \Lambda^n_h,$ & $t\in I$, \\[2mm]
$(\ut_h,\Div~\omega_h)+ (\sigma_{h,t},\omega_h) $ & $= 0,$ & $\forall~\omega_h\in \Lambda^{n-1}_h,$ & $t\in I$, \\[2mm]
$\ut_h(0)$ & $= u_{1,h}$, \\[2mm]
$\sigma(0)$ & $= (\nabla u_0)_h$. \\[2mm]
\end{tabular}
\end{equation}
We now generalize the results of Geveci~\cite{G1988} and others into the language of FEEC.
We first prove a very simple proposition explaining the approximation properties of the $\pi^k_h$ operators in this context.

\begin{proposition}
\label{prop:proj-props}
Choose finite element spaces
\[\Lambda^{n-1}_h = \left\{\begin{array}{c} \poly_{r+1}\Lambda^{n-1}(\Tau) \\ \text{or}\\ \poly_{r+1}^-\Lambda^{n-1}(\Tau)\end{array}\right\},\quad \Lambda^n_h = \poly_{r+1}^-\Lambda^n(\Tau) \;\;\left(= \poly_{r}\Lambda^n(\Tau)\right).
\]
The smoothed projection operators from~(\ref{eq:coch-proj-def})  have the approximation properties
\begin{align}
\vn{\pi^{n-1}_h\omega - \omega}_{L^2\Lambda^{n-1}} & \leq c~h^s\vn{\omega}_{H^s\Lambda^{n-1}},\notag \\
&\begin{array}{ll} 
\text{for $0\leq s\leq r+2$,} & \text{if $\Lambda^{n-1}_h = \poly_{r+1}\Lambda^{n-1}(\Tau)$, or} \\[2mm]
\text{for $0\leq s\leq r+1$,} & \text{if $\Lambda^{n-1}_h = \poly_{r+1}^-\Lambda^{n-1}(\Tau)$,} 
\end{array} \label{eq:omega-approx} \\[2mm]
\vn{\pi^n_h\phi - \phi}_{L^2\Lambda^n} & \leq c~h^s\vn{\phi}_{H^s\Lambda^n},\;\;\text{for $0\leq s\leq r+1$.} \label{eq:ut-approx}
\end{align}
\end{proposition}
\begin{proof}

Estimate (\ref{eq:omega-approx}) follows directly from Theorem~\ref{thm:afw-5pt9}~(\textit{i.}).
Note that Theorem~\ref{thm:afw-5pt9}~(\textit{i.}) is stated for the case $\poly_r\Lambda^{n-1}$ while here we have $\poly_{r+1}\Lambda^{n-1}$, thereby allowing for the higher bound on $s$ in this case.
Finally, since $\pi^n_h\ut-\ut = \p_t(\pi^n_h u-u)$, Theorem~\ref{thm:afw-5pt9}~(\textit{i.}) also implies (\ref{eq:ut-approx}).
\end{proof}

\begin{theorem}
\label{thm:boch-est-hyp}
Fix $\Omega\subset\R^n$ and fix $I:=[0,T]$.  
Choose finite element spaces
\[\Lambda^{n-1}_h = \left\{\begin{array}{c} \poly_{r+1}\Lambda^{n-1}(\Tau) \\ \text{or}\\ \poly_{r+1}^-\Lambda^{n-1}(\Tau)\end{array}\right\},\quad \Lambda^n_h = \poly_{r+1}^-\Lambda^n(\Tau) \;\;\left(= \poly_{r}\Lambda^n(\Tau)\right).
\] 
Then for $(\ut_h,\sigma_h)$ the solution to (\ref{eq:hyp-semidisc}), the following error estimate holds:
\begin{equation}
\label{eq:hyp-err-est-new}
\vn{\ut_h-\ut}_{L^2(I,L^2\Lambda^n)}+ \vn{\sigma_h-\sigma}_{L^2(I,L^2\Lambda^{n-1})} \leq c \left(\sqrt T E_1 + h^s\left(\sqrt TE_2 + E_3 \right)\right),
\end{equation}
where 
\begin{align*}
E_1 & = \vn{u_1-u_{1,h}}_{L^2}+ \\
 & \qquad\vn{(\nabla u_0)-(\nabla u_0)_h}_{L^2} & \text{(error due to discretization of initial data)} \\
E_2 & = \vn{u_1}_{H^s}+\vn{\nabla u_0}_{H^s} & \text{(regularity of initial data)} \\
E_3 & = \vn{u_t}_{L^2(I,H^s)}+\vn{\sigma}_{L^2(I,H^s)} & \text{(regularity of continuous solution to (\ref{eq:boch-hyper-s}))}
\end{align*}
\end{theorem}
\begin{remark}
\em
This theorem strengthens and generalizes the result by Geveci~\cite{G1988} for $n=2$ where $L^2$ projection is used instead of the smoothed projection operators $\pi^k_h$.  
An article by Makridakis~\cite{M1992} extended Geveci's results to $n=3$ in the context of linear elastodynamics, however both papers had to assume the existence of finite element spaces and projections to them with certain properties.
Our result here makes clear what these spaces and projections should be in the unified language of FEEC.
Moreover, the fact that the $\pi^k_h$ operators are not the $L^2$ projection and hence not self-adjoint requires a revised proof technique that ultimately allows the removal of the error term $\vn{u_{tt}}_{L^2(I,H^s)} +\vn{\sigma_t}_{L^2(I,H^s)}$ appearing in prior error bounds.
Recent work by Falk and Winther~\cite{FW2014,FW2015} presents local $L^2$-stable cochain projectors that could potentially be used in place of the $\pi^k_h$ operators.
\end{remark}
\begin{proof}
Define $\Psi:=\Lambda^n\times \Lambda^{n-1}$ with finite dimensional subspace $\Psi_h:=\Lambda^n_h\times \Lambda^{n-1}_h$.
Denote the components of an element $\psi_i\in\Psi$ by $\{\phi_i,\omega_i\}$.  The $L^2$ inner product and norm on $\Psi$ are
\[(\psi_1,\psi_2)_\Psi:= (\phi_1,\phi_2)_{L^2}+(\omega_1,\omega_2)_{L^2}\quad\text{and}\quad\vn{\psi}_\Psi:=\sqrt{(\psi,\psi)_\Psi}.\]
Define a skew-symmetric bilinear form $a:\Psi\times\Psi\raw\R$ by
\[a(\psi_1,\psi_2) := -(\Div~\omega_1,\phi_2)_{L^2}+(\phi_1,\Div~\omega_2 )_{L^2}.\]
Let $\xi:=(\ut,\sigma)\in\Psi$ be the solution to (\ref{eq:boch-hyper-s}) and let $\psi:=(\phi,\omega) \in\Psi$ be arbitrary.
Then adding the equations of (\ref{eq:boch-hyper-s}) yields
\begin{equation}
\label{eq:hyp-single}
(\xi_t,\psi)_\Psi+a(\xi,\psi)=(f,\phi)_{L^2}\qquad\forall~\psi\in\Psi.
\end{equation}
Similarly, from (\ref{eq:hyp-semidisc}) we get
\begin{equation}
\label{eq:hyp-single-disc}
(\xi_{h,t},\psi_h)_\Psi+a(\xi_h,\psi_h)=(f,\phi_h)_{L^2}\qquad\forall~\psi_h\in\Psi_h.
\end{equation}
Define a projection operator $\pi_h:\Psi\raw\Psi_h$ using the bounded cochain projections from (\ref{eq:coch-proj-def}) via $\pi_h\psi:=\left\{\pi^n_h\phi,\pi^{n-1}_h\omega\right\}$.
Since $\pi_h$ only affects the spatial variables, it commutes with the time derivative operator, i.e.\
$(\p_t\pi_h\xi,\psi_h)_{\Psi}=(\pi_h\p_t\xi,\psi_h)_{\Psi}.$
Using this and (\ref{eq:hyp-single}), and letting $\id$ denote the identity operator, we derive
\begin{align}
a(\pi_h\xi,\psi_h) & = a(\xi,\psi_h) + a((\pi_h-\id)\xi,\psi_h) \notag\\
(\p_t\pi_h\xi,\psi_h)_{\Psi} + a(\pi_h\xi,\psi_h) &= a(\xi,\psi_h) + (\pi_h\p_t\xi,\psi_h)_{\Psi} + a((\pi_h-\id)\xi,\psi_h)\notag \\
(\p_t\pi_h\xi,\psi_h)_{\Psi}+a(\pi_h\xi,\psi_h) &  = (f,\phi_h)_{L^2} +\left((\pi_h-\id)\p_t\xi,\psi_h\right)_{\Psi} + a((\pi_h-\id)\xi,\psi_h),\label{eq:hyp-sing-mod}
\end{align}
which holds for all $\psi_h\in\Psi_h\subset\Psi$. 
Now define the error function
\[\veps_h(t) := \pi_h\xi(t)-\xi_h(t).\]
The derivation of a good bound for $\vn{\veps_h(t)}_\Psi$ constitutes the bulk of the remainder of the proof.
Subtracting (\ref{eq:hyp-single-disc}) from (\ref{eq:hyp-sing-mod}) yields
\begin{equation}
\label{eq:geveci-gen}
(\p_t \veps_h(t),\psi_h)_\Psi+a(\veps_h(t),\psi_h)=((\pi_h-\id)\p_t\xi(t),\psi_h)_\Psi + a((\pi_h-\id)\xi(t),\psi_h),\quad\forall~\psi_h\in\Psi_h.
\end{equation}
Define the skew-adjoint linear operator $L_h:\Psi_h\raw\Psi_h$ by
\[\left(L_h\psi_1,\psi_2\right)_{\Psi} := a(\psi_1,\psi_2)\qquad \forall~\psi_1, \psi_2\in\Psi_h.\]
We can thus re-write~(\ref{eq:geveci-gen}) as an equation of functionals on $\Psi_h$:
\begin{equation}
\label{eq:gev-1p19}
\p_t \veps_h(t)+L_h\veps_h(t) = (\pi_h-\id)\p_t\xi(t) + L_h(\pi_h-\id)\xi(t)
\end{equation}
To ease notation, set $Q(t):=(\pi_h-\id)\p_t\xi(t)$ and $R(t):=(\pi_h-\id)\xi(t)$, yielding
\begin{equation}
\label{eq:new-func-eqn}
\p_t \veps_h(t)+L_h\veps_h(t) = Q(t)+L_h R(t)
\end{equation}
We will use some basic results from the theory of semigroups of linear operators as can be found, for instance, in~\cite{P1983}.
For any fixed $\tau\in\R$, the product rule in this context yields
\begin{align}
\p_t(e^{(t-\tau)L_h}(\veps_h(t)-R(t)) 
 &= L_h e^{(t-\tau)L_h}(\veps_h(t)-R(t))+ e^{(t-\tau)L_h} \p_t(\veps_h(t)-R(t)) \notag\\
 &= e^{(t-\tau)L_h}\left(\p_t\veps_h(t)+L_h\veps_h(t) - (L_h+\p_t) R(t) \right) \notag
 \end{align}
Note that we used the fact that $e^{(t-\tau)L_h}$ commutes with $L_h$, a standard result~\cite[Corollary 1.4]{P1983}.
Swapping the roles of $t$ and $\tau$, we re-write the above as
\begin{equation}
\label{eq:semigp-eqn}
\p_\tau(e^{-(t-\tau)L_h}(\veps_h(\tau)-R(\tau))) =  e^{-(t-\tau)L_h}\left(\p_\tau\veps_h(\tau)+L_h\veps_h(\tau) - (L_h+\p_\tau) R(\tau)\right)
\end{equation}
Now we integrate in such a way that (\ref{eq:new-func-eqn}) and (\ref{eq:semigp-eqn}) will give us an expression for $\veps_h(t)$.
First observe that $\p_\tau R(\tau)=Q(\tau)$ since $\p_\tau$ commutes with $\pi_h$ and $\id$.
Thus,
\begin{align*}
0 
& = \int_0^t e^{-(t-\tau)L_h}\left(Q(\tau) - \p_\tau R(\tau)\right)~d\tau \\
&= \int_0^t e^{-(t-\tau)L_h}\left(Q(\tau)+L_hR(\tau) - (L_h+\p_\tau) R(\tau)\right)~d\tau  & \text{by $\pm L_hR(\tau)$}\\
&= \int_0^t e^{-(t-\tau)L_h} \left(\p_\tau \veps_h(\tau)+L_h\veps_h(\tau) - (L_h+\p_\tau) R(\tau)\right)~d\tau  & \text{by (\ref{eq:new-func-eqn})} \\
& = \int_0^t \p_\tau(e^{-(t-\tau)L_h}(\veps_h(\tau)-R(\tau)) ~d\tau & \text{by (\ref{eq:semigp-eqn})} \\
& = \veps_h(t) -R(t) -e^{-t L_h}\veps_h(0) +R(0)& \text{fund.\ thm.\ calculus}
\end{align*}
Rewriting the above chain of equalities, we see that
\begin{equation}
\label{eq:new-veps-exp}
\veps_h(t) = e^{-t L_h}\veps_h(0) + R(t)-R(0)
\end{equation}
Observe that $e^{-tL_h}$ is unitary meaning it preserves $\Psi$-norm, i.e.\ $\vn{e^{-t L_h}\psi}_\Psi = \vn{\psi}_\Psi$ for all $\psi\in\Psi$.  
This follows from the fact that $L_h$ is a real, skew self-adjoint operator, meaning $iL_h$ is self-adjoint, which is equivalent to saying $e^{-tL_h}$ is unitary~\cite[Theorem 10.8]{P1983}.
Thus, taking the $\vn{\cdot}_\Psi$ norm of (\ref{eq:new-veps-exp}), the triangle inequality gives
\begin{align}
\vn{\veps_h(t)}_\Psi 
& \leq \vn{\veps_h(0)}_\Psi + \vn{R(t)}_\Psi +\vn{R(0)}_\Psi \notag \\
& = \vn{\veps_h(0)}_\Psi + \vn{(\pi_h-\id)\xi(t)}_\Psi + \vn{(\pi_h-\id)\xi(0)}_\Psi,\label{eq:new-veps-bd} 
\end{align}
Unpacking the notation lets us characterize this bound in terms of the errors defined in the theorem statement.
Recall that $u_0$ and $u_1$ are given initial data functions and should not be confused with $u_h$ or $u_t$.
We will use $f\lesssim g$ to mean $f\leq cg$ where $c$ is some constant independent of $h$ and $T$.
In $\vn{\cdot}_\Psi$ norm, we then have
\begin{equation}
\label{eq:veps-gen-bd}
\vn{\veps_h(t)} \lesssim \vn{\pi^n_hu_1-u_{1,h}}+\vn{\pi^{n-1}_h(\nabla u_0)-(\nabla u_0)_h} + \vn{(\pi_h-\id)\xi(t)} + \vn{(\pi_h-\id)\xi(0)}  
\end{equation}
To bound the first term on the right, use (\ref{eq:ut-approx}) from Proposition~\ref{prop:proj-props} to get
\[\vn{\pi^n_hu_1-u_{1,h}} \leq \vn{\pi^n_hu_1-u_1}+\vn{u_1-u_{1,h}}\lesssim h^s\vn{u_1}_{H^s}+\vn{u_1-u_{1,h}}\]
Using (\ref{eq:omega-approx}) likewise for the second term, we have
\begin{equation}
\label{eq:init-data-bounds}
\vn{\pi^n_hu_1-u_{1,h}} + \vn{\pi^{n-1}_h(\nabla u_0)-(\nabla u_0)_h} \lesssim E_1 + h^sE_2
\end{equation}
Also by (\ref{eq:omega-approx}) and (\ref{eq:ut-approx}), we have the bounds
\begin{align}
\vn{(\pi_h-\id)\xi(t)} & \lesssim h^s(\vn{\sigma(t)}_{H^s\Lambda^{n-1}}+\vn{u_t(t)}_{H^s\Lambda^{n}}) \label{eq:cnts-bd-t} \\
\vn{(\pi_h-\id)\xi(0)} & \lesssim h^s(\vn{\nabla u_0}_{H^s\Lambda^{n-1}}+\vn{u_1}_{H^s\Lambda^{n}}) \label{eq:cnts-bd-zero}
\end{align}
Using (\ref{eq:veps-gen-bd}) in conjunction with (\ref{eq:init-data-bounds}), (\ref{eq:cnts-bd-t}), and (\ref{eq:cnts-bd-zero}), we derive
\begin{align}
\int_0^T\vn{\veps_h(t)}^2dt 
 & \lesssim T E_1^2 + h^{2s}TE_2^2 + h^{2s}\int_0^T \vn{\sigma(t)}_{H^s}^2+\vn{u_t(t)}_{H^s}^2 + E_2^2~dt \notag \\
 & \lesssim T E_1^2 + h^{2s}TE_2^2 + h^{2s}E_3 +h^{2s}TE_2^2  \label{eq:vepsh-bound}
\end{align}
We now start building up the main result.
\begin{samepage}
\begin{flushleft}
$\ds\left(\vn{\ut_h-\ut}_{L^2(I,L^2\Lambda^n)}+ \vn{\sigma_h-\sigma}_{L^2(I,L^2\Lambda^{n-1})}\right)^2$
\end{flushleft}
\vspace{-.2cm}
\begin{align}
 &\lesssim  \int_0^T \left( \vn{\ut_h-\pi^n_h\ut}+ \vn{\sigma_h-\pi^{n-1}_h\sigma}\right)^2 + \vn{\pi^n_h\ut-\ut}^2+ \vn{\pi^{n-1}_h\sigma-\sigma}^2 dt \notag\\
 & = \int_0^T\vn{\veps_h(t)}^2dt + \vn{\pi^n_h\ut-\ut}_{L^2(I,L^2\Lambda^n)}^2+ \vn{\pi^{n-1}_h\sigma-\sigma}_{L^2(I,L^2\Lambda^{n-1})}^2  dt \notag\\
 & \lesssim \int_0^T\vn{\veps_h(t)}^2 + h^{2s}\vn{\ut(t)}_{H^s}^2 + h^{2s}\vn{\sigma(t)}_{H^s}^2 dt \notag \\
 & \lesssim \int_0^T\vn{\veps_h(t)}^2 dt + h^{2s} E_3^2.  \label{eq:ut-plus-s-err}
\end{align}
\end{samepage}
Combining (\ref{eq:ut-plus-s-err}) and (\ref{eq:vepsh-bound}) yields
\begin{align*}
\left(\vn{\ut_h-\ut}_{L^2(I,L^2\Lambda^n)}+ \vn{\sigma_h-\sigma}_{L^2(I,L^2\Lambda^{n-1})}\right)^2
& \lesssim \left(TE_1^2 + 2h^{2s}TE_2^2 + 2h^{2s}E_3^2 \right) \\
& \lesssim \left(\sqrt TE_1 +h^{s}\left(\sqrt TE_2 + E_3 \right)\right)^2
\end{align*}
Taking the square root of both sides completes the proof. 
\end{proof}

\section{Semi-linear Evolution Problems}
\label{sec:semilin}

We now show how the techniques developed above can be extended to certain types of non-linear evolution problems.  
Consider the \textbf{semi-linear heat equation}: Find $u(x,t)$ such that
\begin{equation}
\label{eq:semilin-cnts}
\begin{tabular}{rlll}
$u_t - \Delta u +F(u) = f$ & in $\Omega$, & for $t>0$ \\
$u = 0$ & on $\p\Omega$, & for $t>0$ & with $u(\cdot,0)=g$ in $\Omega$,
\end{tabular}
\end{equation}
where $F$ is some non-linear operator on $L^2(\Omega)$.   
The existence and uniqueness of solutions to instances of this problem have been studied extensively \cite{G1986,H1981,RT2002,WZ2009} as have finite element methods for the approximation of its solution~
\cite{CH2003,LS2007,S2006,T2006,TW1975}. 

We define the \textbf{semi-linear mixed weak form parabolic} problem: Given $f$ and $g$, find~$(u,\sigma): I \raw H\Lambda^n \times H\Lambda^{n-1}$ such that
\begin{equation}
\label{eq:boch-semilin}
\begin{tabular}{rllll}
$(u_{t},\phi) - (\Div~\sigma,\phi) + (F(u),\phi)$ & $= (f,\phi),$ & $\forall~\phi\in H\Lambda^n,$ & $t\in I$, \\[2mm]
$ (u,\Div~\omega)+ (\sigma,\omega)$ & $= 0,$ & $\forall~\omega\in H\Lambda^{n-1},$ & $t\in I$,  \\[2mm]
$u(0)$ & $= g$.
\end{tabular}
\end{equation} 
The \textbf{semi-linear semi-discrete parabolic problem} is thus: Find $(u_h,\sigma_h): I\raw \Lambda^n_h\times \Lambda^{n-1}_h$ such that
\begin{equation}
\label{eq:par-sl-sd}
\begin{tabular}{rllll}
$(u_{h,t},\phi_h) - (\Div~\sigma_h,\phi_h)+(F(u_h),\phi_h)$ & $= (f,\phi_h),$ & $\forall~\phi_h\in \Lambda^n_h,$ & $t\in I$, \\[2mm]
$(u_h,\Div~\omega_h)+ (\sigma_h,\omega_h) $ & $= 0,$ & $\forall~\omega_h\in \Lambda^{n-1}_h,$ & $t\in I$,\\[2mm]
$u_h(0)$ & $=g_h$,
\end{tabular}
\end{equation}
where $g_h\in \Lambda^n_h$ is an approximation of $g$.  Analogously to the linear case, for any $t_0\in I$, define the \textbf{time-ignorant linear discrete elliptic problem}: find $(\tilde u_h,\tilde\sigma_h)\in \Lambda^n_h\times\Lambda^{n-1}_h$ such that
\begin{equation}
\label{eq:sl-par-to-ellip}
\begin{tabular}{rll}
$(\Div~\tilde\sigma_h,\phi_h) - (\Delta u(t_0),\phi_h)$ & $= 0,$ & $\forall~\phi_h\in \Lambda^n_h,$   \\[2mm] 
$(\tilde\sigma_h,\omega_h) + (\tilde u_h,\Div~\omega_h)$ & $= 0,$ & $\forall~\omega_h\in \Lambda^{n-1}_h,$  \\[2mm] 
$\tilde u_h(0)$ & $= g_h$,
\end{tabular}
\end{equation}
where now $u$ is the solution to the continuous semi-linear problem~(\ref{eq:boch-semilin}).  
Recall, as in (\ref{eq:par-to-ellip}), that no continuity of $\tilde\sigma_h$ or $\tilde u_h$ with respect to $t$ is required.
Similarly, define
\begin{align*}
\rho(t) & := \tilde u_h(t) - u(t),  \\
\theta(t) & := u_h(t) -\tilde u_h(t),  \\
\veps(t) & := \sigma_h(t) -\tilde\sigma_h(t),
\end{align*}
where $u$, $\sigma$ and their discrete counterparts are now solutions to the corresponding semi-linear problems.  
We also introduce the error term
\begin{align*}
\eta(t) & := F(u(t)) - F(u_h(t)).
\end{align*}
Here, 
we assume that the nonlinear function $F: \R \to \R$ satisfies the following (local) Lipschitz condition:
\begin{equation}
\label{eq:F-Lip-assump}
 (F(u)-F(u_{h}), v)\leq K\vn{u-u_{h}}_{L^{2}} \|v\|,\quad\forall~v\in L^2(\Omega),
\end{equation}
where $u$ and $u_{h}$ are the solutions to \eqref{eq:boch-semilin} and \eqref{eq:par-sl-sd} respectively. 
We remark that this type of local Lipschitz condition is slightly different from those assumed by Holst and Stern in~\cite{HS2010b,BHSZ2011}.  
Here, we use the $L^{2}$ norm on the right sides of the inequalities. 
By the Poincar\'e inequality, \eqref{eq:F-Lip-assump} implies the Lipschitz assumption used in~\cite{HS2010b,BHSZ2011}.  
The key here is that we only require the assumption \eqref{eq:F-Lip-assump} to hold for the solutions $u$ and $u_{h}$.  
For example, suppose that the solutions $u$ and $u_{h}$ satisfy some \emph{a priori} $L^{\infty}$ bounds i.e.\ there exist some constants $\alpha$ and $\beta$ such that $\alpha \le u, u_{h} \le \beta$.
Then the local Lipschitz condition is satisfied if there is some constant $K >0$ such that $|F'(\xi)| \le K$ for all $\xi \in [\alpha, \beta].$  

We have the follow lemma, which reduces to Lemma~\ref{lem:err-eqns-lin} in the linear case.
\begin{lemma}
\label{lem:err-eqns-sl}
The semi-linear error functions satisfy the semi-discrete formulation, i.e.
\begin{equation}
\label{eq:par-error-sl}
\begin{tabular}{rllll}
$(\theta_{t},\phi_h) - (\text{\em div}~\veps,\phi_h) $ & $= (\eta,\phi_h)-(\rho_t,\phi_h),$ & $\forall~\phi_h\in \Lambda^n_h,$ & $t\in I$, \\[2mm]
$(\theta,\text{\em div}~\omega_h)+(\veps,\omega_h) $ & $= 0,$ & $\forall~\omega_h\in \Lambda^{n-1}_h,$ & $t\in I$.
\end{tabular}
\end{equation}
\end{lemma}
\begin{proof}
The second equation is immediate from the second equations in (\ref{eq:par-sl-sd}) and (\ref{eq:sl-par-to-ellip}).  
The first equation is derived as follows.
\begin{align*}
(u_{t},\phi_h) - (\Delta u,\phi_h) + (F(u),\phi_h) & = (f,\phi_h) &  \text{by  (\ref{eq:semilin-cnts}),} \\
(u_{t},\phi_h) - (\Delta u,\phi_h) + (F(u),\phi_h) & = (u_{h,t},\phi_h) - (\Div~\sigma_h,\phi_h)+(F(u_h),\phi_h) &  \text{by  (\ref{eq:par-sl-sd}),} \\
(u_{t},\phi_h) - (\Div~\tilde\sigma_h,\phi_h) + (F(u),\phi_h) & = (u_{h,t},\phi_h) - (\Div~\sigma_h,\phi_h)+(F(u_h),\phi_h) &  \text{by  (\ref{eq:sl-par-to-ellip})},
\end{align*}
from which the result follows by adding $(\pm~\tilde u_{h,t},\phi_h)$ and collecting terms.
\end{proof}

\begin{theorem}
\label{thm:semilin-est}
Fix $\Omega\subset\R^n$ and fix $I:=[0,T]$. 
Suppose $(u,\sigma)$ is the solution to (\ref{eq:boch-semilin}) such that the regularity estimate
\begin{equation}
\label{eq:reg-est-sl}
\vn{u(t)}_{H^{s+2}}+\vn{\sigma(t)}_{H^{s+1}}+\vn{d\sigma(t)}_{H^{s}}\leq c\vn{\Delta u(t)}_{H^{s}}
\end{equation}
holds for $0\leq s\leq s_{\max}$ and $t\in I$.  
Assume that the operator $F$ satisfies the Lipschitz assumption (\ref{eq:F-Lip-assump}).
Choose finite element spaces 
\[\Lambda^{n-1}_h = \left\{\begin{array}{c} \poly_{r+1}\Lambda^{n-1}(\Tau) \\ \text{or}\\ \poly_{r+1}^-\Lambda^{n-1}(\Tau)\end{array}\right\},\quad \Lambda^n_h = \poly_{r+1}^-\Lambda^n(\Tau) \;\;\left(= \poly_{r}\Lambda^n(\Tau)\right)
\]
Then for $0\leq s\leq s_{\max}$, $g_h$ defined by (\ref{eq:gh-def}), and $(u_h,\sigma_h)$ the solution to (\ref{eq:par-sl-sd}), the following error estimates hold:
\begin{align}
\vn{u_h-u}_{L^2(I,L^2\Lambda^n)} &\lesssim\left\{\begin{array}{ll} \ds  h\left(\vn{\Delta u}_{L^2(I,L^2)}+\vn{\Delta u_t}_{L^2(I,L^2)}\right) & \text{if $r=0$} \\[4mm] \ds h^{1+s}\left(\vn{\Delta u}_{L^2(I,H^s)}+ \vn{\Delta u_t}_{L^2(I,H^s)}\right) & \text{for $r>0$,} \\
 & \text{~~if $s\leq r-1$}\end{array}\right.\label{eq:semilin-est-u}
\end{align}
\begin{align}
\vn{\sigma_h-\sigma}_{L^2(I,L^2\Lambda^{n-1})} & \lesssim \left\{\begin{array}{l} \ds h\left(\vn{\Delta u}_{L^2(I,H^s)}+\vn{\Delta u_t}_{L^2(I,L^2)}\right) \\[3mm] \qquad \text{if~$r=0$, $s=0$, $\Lambda^{n-1}_h=\poly_1^-\Lambda^{n-1}(\Tau)$} \\[3mm]
\ds h^{1+s}\vn{\Delta u}_{L^2(I,H^s)}+h\vn{\Delta u_t}_{L^2(I,L^2)} \\[3mm] \qquad\text{if~$r=0$, $s\leq 1$, $\Lambda^{n-1}_h=\poly_1\Lambda^{n-1}(\Tau)$} \\[3mm]
\ds h^{1+s}\vn{\Delta u}_{L^2(I,H^s)}+h^{(3/2)+s}\vn{\Delta u_t}_{L^2(I,H^s)} \\[3mm] \qquad\text{for $r>0$, if $s\leq r-1$}
\end{array}\right. \label{eq:semilin-est-s}
\end{align}
\begin{align}
\vn{\textnormal{div}(\sigma_h-\sigma)}_{L^2(I,L^2\Lambda^n)} & \lesssim\left\{\begin{array}{l} \ds h^s\vn{\Delta u}_{L^2(I,H^s)}+h\vn{\Delta u_t}_{L^2(I,L^2)} \\[3mm] \qquad \text{if~$r=0$, $s\leq 1$}\\[3mm]
\ds h^s\vn{\Delta u}_{L^2(I,H^s)}+h^{2+s}\vn{\Delta u_t}_{L^2(I,H^s)} \\[3mm] \qquad\text{for $r>0$, if $s\leq r-1$}
\end{array}\right..\label{eq:semilin-est-ds}
\end{align}
\end{theorem}
\begin{proof}
The proof is very similar to that of Theorem~\ref{thm:boch-est-par}.  Equation~(\ref{eq:sl-par-to-ellip}) is the $k=n$ case of the discrete mixed variational problem examined by Holst and Stern in~\cite[Equation (9)]{HS2010b}.  Therefore, we can use the same type of triangle inequality from (\ref{eq:u-err-tri-ineq}) to recover the estimates.  By~\cite[Theorem 4.2]{HS2010b}, we have the estimates
\begin{align}
\|\rho(t) \| & \leq\left\{\begin{array}{ll} ch \|\Delta u(t) \| & \text{if $r=0$} \\[2mm] ch^{2+s}\vn{\Delta u(t)}_{H^s} & \text{if $s\leq r-1$, for $r>0$}\end{array}\right.\label{eq:rho-est-sl} \\
\|\rho_{t}(t)\| & \leq\left\{\begin{array}{ll} ch \|\Delta u_{t}(t)\| & \text{if $r=0$} \\[2mm] ch^{2+s}\vn{\Delta u_{t}(t)}_{H^s} & \text{if $s\leq r-1$, for $r>0$}\end{array}\right.\label{eq:rhot-est-sl} \\
\|\sigma(t)-\tilde\sigma_h(t)\| & \leq ch^{1+s}\vn{\Delta u(t)}_{H^s},~\text{if}~\left\{\begin{array}{ll} s\leq r+1, & \Lambda^{n-1}_h=\poly_{r+1}\Lambda^{n-1}(\Tau) \\[2mm] 
s\leq r, & \Lambda^{n-1}_h=\poly_{r+1}^-\Lambda^{n-1}(\Tau) \end{array}\right.\label{eq:sigma-est-sl} \\
\|\Div(\sigma(t)-\tilde\sigma_h(t))\| & \leq ch^{s}\vn{\Delta u(t)}_{H^s},~\text{if~$s\leq r+1$}.
\label{eq:dsigma-est-sl}
\end{align}
An explanation of how these estimates are derived from the results of~\cite{HS2010b} is given in Appendix~\ref{app:semilin-est}.  
Note that these estimates are exactly the same as the corresponding estimates (\ref{eq:rho-est}), (\ref{eq:AFWsigma-est}) and (\ref{eq:AFWdsigma-est}) from the linear case. 

To bound $\vn{\theta}$, we adapt the proof of Theorem~\ref{thm:boch-est-par} to account for the non-linearity.
Set $\phi_h := \theta$ and $\omega_h := \veps$ in (\ref{eq:par-error-sl}).
Adding the equations yields
\begin{equation}
\label{eq:norms-to-ip-sl}
\frac 12 \frac{d}{dt}\|\theta\|^2 + \|\veps\|^2 + (F(u_h) - F(u), \theta) = -(\rho_t,\theta),\;\; t\in I
\end{equation}
By the Lipschitz condition \eqref{eq:F-Lip-assump}, we obtain: 
\begin{align*}
	\frac 12 \frac{d}{dt}\|\theta\|^2  &\le |(F(u) - F(u_{h}), \theta)|  + |(\rho_t,\theta)|\\
	&\le K\|u-u_{h}\| \|\theta\|+ \|\rho_{t}\| \|\theta\| \\
	&\le K\|\rho\|\|\theta\| + \|\rho_{t}\| \|\theta\| + K \|\theta\|^{2}\\
	&\le C\|\theta\|^{2} + C(\|\rho\|^{2} + \|\rho_{t}\|^{2}) .
\end{align*}
Notice that $\theta(0) =0$ and an application of Gronwall's inequality implies that 
\begin{align*}
	\|\theta(t)\|^{2} &\le 2C e^{2Ct} \int_{0}^{t} e^{-2Cs} (\|\rho\|^{2} + \|\rho_{s}\|^{2}) ds\\
	&\le 2C e^{2Ct} \left(\|\rho\|_{L^{2}(I, L^{2})}^{2} + \|\rho_{t}\|_{L^{2}(I, L^{2})}^{2} \right).
\end{align*}
To get an estimate on $\|\theta\|_{L^{2}(I, L^{2})}$, we integrate by $t$ on both  sides to obtain 
\begin{equation}
\label{eqn:theta}
	\|\theta\|_{L^{2}(I, L^{2})} \lesssim \sqrt{T} (\|\rho\|_{L^{2}(I,L^{2})} + \|\rho_{t}\|_{L^{2}(I,L^{2})}).
\end{equation}
The estimate \eqref{eq:semilin-est-u} then follows by \eqref{eq:rho-est-sl}, \eqref{eq:rhot-est-sl}, \eqref{eqn:theta} and a triangle inequality.

To obtain the error bound \eqref{eq:semilin-est-s} for the approximation of $\sigma$,  we use the same technique as in Theorem~\ref{thm:boch-est-par}. To bound $\vn{\eps(t)}$, we differentiate the second equation of \eqref{eq:par-error-sl} with respect to $t$ and set $\phi_h := 2\theta_t$, $\omega_h := 2\veps$, yielding
\begin{equation*}
\begin{tabular}{rllll}
$(\theta_{t},2\theta_t) - (\Div~\veps,2\theta_t)$ & $= (\eta, 2\theta_{t})-(\rho_t,2\theta_t),$\\[2mm]
$(\veps_t,2\veps) + (\theta_t,2\Div~\veps)$ & $= 0.$ 
\end{tabular}
\end{equation*}
Adding the equations and using the Lipschitz assumption \eqref{eq:F-Lip-assump}, we have the bound
\begin{align*}
	\frac {d}{dt} \vn{\veps}^2+ 2\vn{\theta_t}^2 &= 2(\eta, \theta_{t})- 2(\rho_t,\theta_t) \\
	&\leq 2K\|u-u_{h}\|\|\theta_{t}\| + 2\vn{\rho_t}\vn{\theta_t}
\end{align*}	
Now write $\vn{u-u_h}=\vn{\theta+\rho}$, apply the triangle inequality, and then apply Young's inequality to get a constant $C$ such that
\[\frac {d}{dt} \vn{\veps}^2+ 2\vn{\theta_t}^2 \leq C\left(\|\rho\|^{2} + \|\theta\|^{2} + \|\rho_{t}\|^{2} \right) + \|\theta_{t}\|^{2}.\]
Canceling one of the $\|\theta_{t}\|^{2}$ on both sides, and by the estimates \eqref{eq:rho-est-sl}-\eqref{eq:rhot-est-sl} and \eqref{eqn:theta}, we obtain
\begin{align*}
\vn{\veps(t)}^2 & = \int_0^t\frac {d}{ds}\vn{\veps(s)}^2ds  \leq \int_0^t \frac {d}{ds}\vn{\veps(s)}^2 +\vn{\theta_t}^2 ds \\
  &\leq C\int_0^t \|\rho\|^{2} + \|\theta\|^{2} + \|\rho_{t}\|^{2} ds.
\end{align*}
Now use the estimates (\ref{eq:rho-est-sl}) and (\ref{eq:rhot-est-sl}) to derive
\begin{equation}
\label{eq:veps-est-sl}
\vn{\veps(t)} \leq 
\left\{\begin{array}{ll} \ds ch\sqrt{T}\left(\int_0^t \|\Delta u(\ell)\|+ \|\Delta u_t(\ell)\| d\ell\right)^{1/2} & \text{if $r=0$,} \\[4mm] \ds ch^{2+s}\sqrt{T}\left(\int_0^t\vn{\Delta u(\ell)}_{H^s}^2+ \vn{\Delta u_t(\ell)}_{H^s}^2 d\ell\right)^{1/2} & \text{for $r>0$, if $s\leq r-1$.}\end{array}\right.
\end{equation}
The estimate \eqref{eq:semilin-est-s} then follows from triangle inequality, \eqref{eq:sigma-est-sl}, and \eqref{eq:veps-est-sl}.

Finally, to estimate \eqref{eq:semilin-est-ds}, we follow the same techniques. We first bound $\vn{\Div~\veps}$ by setting $w_h:=\veps$ in~(\ref{eq:par-error-sl}) and take the derivative with respect to $t$.  This yields
\[\frac{d}{dt}\vn{\veps}^2+(\theta_t,\Div~\veps)=0.\]
Note that $\Div~\veps\in \Lambda^n_h$ since the discrete spaces are chosen to satisfy the relationship $\Div\Lambda^{n-1}_h\subset\Lambda^n_h$.  Thus, we can set $\phi_h :=\Div~\veps$ in \eqref{eq:par-error-sl} and substitute to get
\begin{equation}
\label{eq:div-eps-est-sl}
\frac{d}{dt}\vn{\veps}^2+\vn{\Div~\veps}^2 = (\rho_t,\Div~\veps) - (\eta,\Div~\veps).
\end{equation}
Integrating both sides over $[0,T]$ and applying the Lipschitz assumption \eqref{eq:F-Lip-assump} and Young's inequality, we have
\begin{equation}
\vn{\veps(T)}^2- \vn{\veps(0)}^2 + \int_0^T\vn{\Div~\veps}^2 dt  \leq C\int_0^T\vn{\rho_t}^2 +\|\rho\|^{2} + \|\theta\|^{2} dt + \frac{1}{2}\int_{0}^{T}\vn{\Div~\veps}^2 dt.
\end{equation}
Note that $\veps(0)=\sigma_h(0)-\tilde\sigma_h(0)=0$ since the elliptic projection is the identity at time $t=0$.
We thus have
\begin{equation}
\int_0^T\vn{\Div~\veps}^2 dt  \leq 2 C\int_0^T\vn{\rho_t}^2 +\|\rho\|^{2} + \|\theta\|^{2} dt.
\end{equation}
The rest of the proof follows similarly.
\end{proof}

We remark that the estimates for the semi-linear case given by Theorem~\ref{thm:semilin-est} are very similar to those for the linear case given by Theorem~\ref{thm:boch-est-par}, although the latter are somewhat sharper.
In the semi-linear case, the constant on the right side of the equation will have a factor of $\sqrt{T}$, while in the linear case the dependence on $\sqrt{T}$ is made explicit.
Further note that the estimate on $\vn{u_h-u}$ in the linear case involves an $L^1(\cdots\!)$ norm on $\Delta u_t$ as part of the bound, which is a stronger result than the bound we obtained in the semi-linear case that involves the $L^2(\cdots\!)$ norm on $\Delta u_t$.

\section{Concluding Remarks}
\label{sec:conc}

  In this article, we have extended the Finite Element Exterior Calculus
  of Arnold, Falk, and Winther~\cite{AFW2010,AFW2006} 
  for linear mixed variational problems to linear and semi-linear parabolic 
  and hyperbolic evolution systems.
  Both the parabolic and hyperbolic cases make strong use of the smoothed 
  projection operators $\pi^k_h$, which are one of the most elaborate and 
  delicate constructions in the FEEC framework.
  In the parabolic case, the use of the $\pi^k_h$ operators was hidden 
  somewhat by the use of elliptic projection error estimates, proofs of which
  rely on properties of these operators.
  In the hyperbolic case, the proof techniques use these properties more 
  explicitly.
  In both cases, the formal treatment and generalization of these operators by 
  Arnold, Falk and Winther can now be seen as a useful tool for the analysis
  of evolution problems as well as elliptic PDE.
  
  We have also seen in this article how the recent generalizations of the FEEC 
  by Holst and Stern~\cite{HS2010a,HS2010b} for semi-linear elliptic PDE can 
  be extended to evolution PDE as well, both parabolic and hyperbolic types.
  We also anticipate that the basic approach to analyzing variational
  crimes in~\cite{HS2010a,HS2010b} for the linear and semilinar elliptic cases
  will also work in the case of evolution problems; we will explore the
  question of variational crimes in a subsequent article, with the target
  being the analysis of surface finite element methods for evolution problems.

\appendix
\begin{appendix}

\section{Explanation of Semi-Linear Error Estimates}
\label{app:semilin-est}
In this appendix, we explain why estimates (\ref{eq:rho-est-sl}), (\ref{eq:sigma-est-sl}), and (\ref{eq:dsigma-est-sl}) follow from~\cite[Theorem 4.2]{HS2010b}.  We will focus just on the $r>0$ case of (\ref{eq:rho-est-sl}) as it requires the sharpening of a special case of an estimate appearing in~\cite[Theorem 4.2]{HS2010b}. The other cases work out along similar lines by a direct application of the Holst and Stern estimates.

First, we recall some notation from~\cite{AFW2010} used in~\cite{HS2010b}.  If $(W,d)$ is a Hilbert complex with associated domain complex $(V,d)$ and parametrized subcomplex family $(V_h,d)$, denote the best approximation in $W$-norm by
\[E(w) = \inf_{v\in V^k_h}\vn{w-v}_W,\quad w\in W^k.\]
The relevant result from Holst and Stern~\cite[Theorem 4.2]{HS2010b} is stated as 
\begin{align}
\vn{u-\tilde u_h}_V+\vn{p-p_h}_W & \leq c(E(u) +E(du) + E(p) \notag \\
 &\quad + \eta[E(\sigma)+E(d\sigma)] + (\delta+\mu)E(d\sigma)+\mu E(P_{\mathfrak B} u)),\label{eqn:HS-u-est}
\end{align}
where $\eta$, $\delta$, and $\mu$ are coefficients defined as the norms of certain abstract operators, $u\in W_k$, and $p$ is a harmonic $k$-form with discrete counterpart $p_h$ introduced to make the abstract Hodge-Laplacian problem well-posed.  

Casting this into the context of the deRham complex, we have
\[(W,d)=(L^2\Lambda,d)\quad\text{and}\quad (V,d)=(H\Lambda,d).\]
Since we are interested here only in the case $k=n$, there are no harmonic $k$-forms so that $p=p_h=0$.  Further, $du=0$ since $d\Lambda^n=0$, whereby $\vn{u-\tilde u_h}_V=\vn{u-\tilde u_h}_W=\eltn{u-\tilde u_h}$.  This eliminates the error terms in $p$ and $du$, giving us the reduced estimate 
\[\eltn{u-\tilde u_h}\leq c(E(u) + \eta[E(\sigma)+E(d\sigma)] + (\delta+\mu)E(d\sigma)+\mu E(P_{\mathfrak B} u)). \]
Crucially, this estimate can be reduced further when $k=n$.  The derivation of (\ref{eqn:HS-u-est}) uses the estimate 
\[\vn{d(u-\tilde u_h)}_W \leq c(E(du)+\eta[E(d\sigma)+E(p)]\]
from~\cite[Theorem 2.11]{AFW2010} which is unnecessary here since the left side is always zero.  Since this is the only part of the derivation that requires the term $\eta E(d\sigma)$, we can drop it, yielding
\begin{equation}
\label{eq:reduced-HS-est}
\eltn{u-\tilde u_h}\leq c(E(u) + \eta E(\sigma) + (\delta+\mu)E(d\sigma)+\mu E(P_{\mathfrak B} u)).
\end{equation}
We now give bounds on each of the terms in~(\ref{eq:reduced-HS-est}).  The coefficients appearing in the abstract estimates can be stated in terms of powers of $h$ in the deRham context.  These appear in~\cite[p. 312]{AFW2010} as 
\[\eta = O(h),\quad \delta = O(h^{\min(2,r+1)}),\quad\text{and~}\mu=O(h^{r+1}).\]
To bound the error terms, Arnold, Falk and Winther define smooth projection operators $\pi^k_h:L^2\Lambda^k(\Omega)\raw\Lambda^k_h$ satisfying optimal convergence rates as stated precisely in~\cite[Theorem 5.9]{AFW2010}.  For instance, if $\Lambda^k_h$ is one of $\poly_{r+1}^-\Lambda^k(\Tau_h)$ or, if $r\geq 1$, $\poly_r\Lambda^k(\Tau_h)$ then 
\[\vn{w-\pi^k_hw}_{L^2\Lambda^k(\Omega)}\leq ch^s\vn{w}_{H^s\Lambda^k(\Omega)},\quad\text{for~} w\in H^s\Lambda^k(\Omega),\quad 0\leq s\leq r+1.\]
These types of results bound $E(w)$ in terms of $\vn{w}_{H^s\Lambda^k}$, which is in turn bounded in terms of $\vn{\Delta u}_{H^s\Lambda^k}$ by the regularity hypothesis (\ref{eq:reg-est-sl}).  Summarizing these results, we have
\begin{align*}
E(u) & \leq ch^{s+2}\vn{\Delta u}_{H^s}\\
E(\sigma) & \leq ch^{s+1}\vn{\Delta u}_{H^s}\\
E(d\sigma) & \leq ch^s\vn{\Delta u}_{H^s}\\
E(P_{\mathfrak B} u) & \leq ch^{s+2}\vn{\Delta u}_{H^s}
\end{align*}
We can now prove~(\ref{eq:rho-est-sl}) by collecting results and applying them to~(\ref{eq:reduced-HS-est}), yielding  
\[\eltn{\rho(t)}\leq c(h^{s+2}+h(h^{s+1})+(h^{\min(2,r+1)}+h^{r+1})h^s+h^{r+1}h^{s+2})\vn{\Delta u(t)}_{H^s}\]
The greatest common factor from the above expression is $h^{s+2}$ hence this is the overall order estimate that can be inferred, as was claimed.
\end{appendix}

\noindent
\textbf{Acknowledgments.}
MH was supported in part by NSF Awards~0715146,
by DOD/DTRA Award HDTRA-09-1-0036, and by NBCR.
AG was supported in part by NSF Award~0715146 and by NBCR.
YZ was supported in part by NSF DMS 1319110. 
The authors would like to thank an anonymous referee for pointing out some important subtle properties of the $\pi^k_h$ operators as well providing various useful suggestions for technical simplifications of the proofs.
\vspace{-5mm}

\bibliographystyle{abbrv}
\bibliography{../bib/mjh,../bib/akg,../bib/HoSt2010b}

\vspace*{0.5cm}

\end{document}